\documentclass{amsart}
\usepackage{amssymb,amsmath,amsfonts,amscd,euscript}

\newcommand{\C}{\mathbb C}

\newcommand{\N}{\mathbb N}
\newcommand{\Z}{\mathbb Z}
\newcommand{\Pro}{\mathbb P}

\newcommand{\T}{\otimes}
\newcommand{\Cn}{\C[e_0,\ldots,e_{n-1}]}
\newcommand{\an}{a_1,\ldots,a_n}

\newcommand{\ld}{\ldots}

\newcommand{\hk}{\hookrightarrow}
\newcommand{\hki}{\hookleftarrow}
\newcommand{\sur}{\twoheadrightarrow}

\newcommand{\Prod}{\prod\limits}

\newcommand{\nc}{\newcommand}

\nc{\pin}{\pi_1,\ld,\pi_n}
\nc{\zn}{z_1,\ld,z_n}
\nc{\slt}{\mathfrak{sl}_2}
\nc{\g}{\mathfrak{g}}
\nc{\Sh}{{\bf Sh}}
\nc{\slth}{\widehat{\slt}}
\nc{\Slth}{\widehat{\Slt}}
\nc{\Tpi}{\pi_1\T\ld\T\pi_n}
\nc{\fpi}{\pi_1*\ld*\pi_n}
\nc{\A}{\mathfrak A}
\nc{\Vg}{\mathfrak V}
\nc{\Ug}{\mathfrak U}
\nc{\Ga}{\mathfrak G}

\nc{\yn}{y_1,\ld,y_n}
\nc{\pz}{\phi_Z}
\nc{\al}{\alpha}
\nc{\be}{\beta}
\nc{\CT}{(\C^2)^{\T n}}
\nc{\Ca}{\C^{a_1} *\ld *\C^{a_n}}
\nc{\veps}{\varepsilon}
\nc{\ik}{i_1,\ld,i_k}
\nc{\ws}{\widetilde\sigma}
\nc{\W}{\widetilde W}

\nc{\Pn}{P_1,\ld,P_n}
\nc{\ov}{\overline}

\nc{\bs}{b_1,\ldots,b_s}

\nc{\qb}[2]{{\genfrac{[}{]}{0pt}{0}{#1}{#2}}_q}
\nc{\gf}[2]{\genfrac{}{}{0pt}{}{#1}{#2}}
\nc{\qf}[1]{({#1})_q!}

\nc{\pa}{\partial}

\nc{\codim}{{\mathop{\rm codim}}}
\nc{\Id}{{\mathop{\rm Id}}}
\nc{\Tr}{{\mathop{\rm Tr}}}
\nc{\ch}{{\mathop{\rm ch}}}
\nc{\V}{{\mathop{\rm V}}}
\nc{\ann}{{\mathop{\rm Ann}}}
\nc{\im}{{\mathop{\rm im}}}
\nc{\const}{{\mathop{\rm const}}}
\nc{\sltc}{\slt\otimes\C[t]}
\nc{\bra}{\langle}
\nc{\ket}{\rangle}
\nc{\Slt}{{\mathop{\rm SL}}_2}
\nc{\G}{{\mathop{\rm G}}}
\nc{\U}{{\mathop{\rm U}}}
\nc{\Ho}{{\mathop{\rm H}}}
\nc{\Vect}{{\mathop{\rm Vect}}}

\nc{\Sltc}{\Slt (\C[t]/t^n)}
\nc{\sh}{{\mathop{\rm sh}}}
\nc{\diag}{{\mathop{\rm diag}}}
\nc{\Pic}{{\mathop{\rm Pic}}}
\nc{\Lg}{{\mathop{\rm L}}}
\nc{\Fl}{{\mathop{\rm Fl}}}

\nc{\shn}{\sh^{(n)}}
\nc{\shin}{\sh_i^{(n-1)}}

\nc{\Ob}{\EuScript {O}}
\nc{\Lie}{\EuScript {L}}
\nc{\E}{\EuScript {E}}
\nc{\F}{\EuScript {F}}

\newtheorem{rem}{Remark}[section]
\newtheorem*{rem*}{Remark}
\newtheorem{opr}{Definition}[section]
\newtheorem{theorem}{Theorem}[section]
\newtheorem*{theorem*}{Theorem}

\newtheorem{lem}{Lemma}[section]
\newtheorem{cor}{Corollary}[section]
\newtheorem{prop}{Proposition}[section]

\begin{document}
\pagestyle{plain}
\markboth{B.Feigin and E.Feigin}
{Schubert varieties and the fusion products}
\title[Schubert varieties and the fusion products]
{Schubert varieties and the fusion products}
\thanks{2000 Mathematics Subject Classification 17B67}
\author{B.Feigin and E.Feigin}
\thanks{BF: Landau Institute for Theoretical Physics, Russia,
Chernogolovka, 142432,\\ 
feigin@mccme.ru}
\thanks{EF: Independent University of Moscow, Russia, Moscow,
Bol'\-shoi Vlas\-'ev\-skii per.,7, and \\
Moscow State University, Russia, Moscow, Leninskie gori, 1,\\
evgfeig@mccme.ru}
\date{}
\begin{abstract}
For each $A\in\N^n$ we define a Schubert variety $\sh_A$ as a closure
of the $\Slt(\C[t])$-orbit in the projectivization of the fusion product
$M^A$.
We clarify the connection of the geometry of the Schubert varieties
with an algebraic structure of $M^A$ as $\sltc$ modules.
In the case, when all the entries of $A$
are different,
$\sh_A$ is smooth projective complex algebraic
variety.
We study its geometric properties: the Lie algebra of the vector
fields, the coordinate ring, the cohomologies of the line bundles.
We also prove that the fusion products can be realized as the dual spaces of
the sections of these bundles. 
\end{abstract}

\maketitle
\section*{Introduction}
The main goal of the paper is the study of the closure of the orbit of
the lowest weight vector in the fusion product of the $\slt$ modules. Let us
recall the main definitions (see \cite{fusion, mi1}).

We start with an abelian situation. Consider an abelian one-dimensional
Lie algebra with a basis $e$. Let $\C^{a_i}\simeq \C[e]/e^{a_i}$,
$i=1,\ld,n$ be its cyclic modules with a cyclic vector $1$.
For $Z=(z_1,\ld,z_n)\in\C^n$ with a distinct entries consider
the tensor product $\bigotimes_{i=1}^n \C^{a_i}(z_i)$ of the evaluation
representations of the Lie algebra $e\T \C[t]$ (recall that
for the representation $V$ of the Lie algebra $\g$ and the complex number $z$
the evaluation module $V(z)$ of the current algebra $\g\T\C[t]$
coincides with $V$ as a vector space and the action is
given by the formula $x_j\cdot v= z^j (x\cdot v)$
for $x_j=x\T t^j, x\in\g$, $v\in V$). Surely,
$\bigotimes_{i=1}^n \C^{a_i}(z_i)$
is cyclic $\C[e_0,\ld,e_{n-1}]$ module with a cyclic vector being the product
of the cyclic vectors of $\C[e]/e^{a_i}$ and thus is isomorphic to the
quotient of the ring $\C[e_0,\ld,e_{n-1}]$ by some ideal $I(Z)$.
It was proved in
\cite{mi1} that the family of the ideals $I(Z)$ can be continuously extended
to any point of the configuration space $\C^n$. Thus, for any $Z\in\C^n$
we have a quotient $M^A(Z)=\C[e_0,\ld,e_{n-1}]/I(Z)$, where
$A=(a_1,\ld,a_n)\in\N^n$.

Let $\A_n$ be an abelian Lie group, generated by
$\exp (t_ie_i), t_i\in\C, 0\le i< n$.
For a vector $v\in V$ denote by $[v]\in\Pro(V)$ the line $\C\cdot v$.
Define a Schubert variety $\sh_A(Z)\hk\Pro (M^A(Z))$ as a closure of the
$\A_n$-orbit of the point $[u]$, where $u$ is a cyclic vector in $M^A(Z)$:
$$\sh_A(Z)=\ov{\A_n\cdot [u]}\hk \Pro (M^A(Z)).$$
Surely, for the $Z$ with pairwise distinct $z_i$
$\sh_A(Z)\simeq\underbrace{\Pro^1\times\ld\times\Pro^1}_n.$
Now, let $i_1+\ld +i_s=n$ and
$$Z=(\underbrace{z_1,\ld,z_1}_{i_1},\underbrace{z_2,\ld,z_2}_{i_2},\ld,
\underbrace{z_s,\ld,z_s}_{i_s}),\qquad z_i\ne z_j.$$
In \cite{mi1} it was proved that
\begin{equation}
\label{cite}
M^A(Z)\simeq M^{(a_1,\ld,a_{i_1})}(z_1,\ld,z_1)\T\ld
\T M^{(a_{n-i_s+1},\ld,a_n)}(z_s,\ld,z_s).
\end{equation}
It was also shown in \cite{mi1} that for any $z\in\C$ $M^A(z,\ld,z)$ is
isomorphic to $M^A(0,\ld,0)$ and the isomorphism is given by the shift
$z\to 0$.
Denote $\sh_A=\sh_A(0,\ld,0)$. From the formula $(\ref{cite})$ we obtain
$$\sh_A(Z)\simeq \sh_{(a_1,\ld,a_{i_1})}\times\ld\times
\sh_{(a_{n-i_s+1},\ld,a_n)}.$$
(Note that the family $\sh_A(Z)$ defines a variety $\Sh_A$, which is fibered
over
$\C^n$ with a fiber over a point $Z$ being $\sh_A(Z)$, but we do not study
$\Sh_A$ in this paper).

Recall that in \cite{mi1} $M^A(0,\ld,0)$ was identified with a fusion
product of the $\slt$ modules. Thus the varieties $\sh_A$ can be defined
in terms of the fusion product. Recall the main points.

Let $e,h,f$ be the standard $\slt$ basis,
$\C^{a_1},\ld,\C^{a_n}$ -- irreducible $\slt$ modules, $z_1,\ld,z_n$ --
a set of the pairwise distinct complex numbers.
Surely,
the tensor product $\bigotimes \C^{a_i}(z_i)$ is a cyclic
$\slt\T\C[t]$ module with a cyclic vector $\ov{v_A}$ being the tensor
product
of the lowest weight vectors (with respect to the $h$ grading). Thus we
have a filtration $F_s$ on $\bigotimes \C^{a_i}(z_i)$ (recall that
$x_j=x\T t^j$)
$$F_s=\bra x^{(1)}_{i_1}\ld x^{(k)}_{i_k}\ov{v_A}:\ x^{(j)}\in\slt,
i_1+\ld + i_k\le s\ket.$$
The adjoint
graded module is denoted by $\C^{a_1}*\ld *\C^{a_n}$ or simply $M^A$ for
$A=(a_1,\ld,a_n)$. Let $v_A\in M^A$ be the image of $\ov{v_A}$. Surely,
$$M^A=U(\slt\T \C[t])\cdot v_A=U(e\T \C[t])\cdot v_A$$
(in fact, the operators from $\slt\T t^i, i\ge n$ vanish on $M^A$).

Consider a group $\Slt (\C[t])$ acting on $M^A$, and thus on its
projectivization $\Pro (M^A)$. The
Schubert variety can be defined as follows: it is the closure of the
$\Slt (\C[t])$-orbit of the point $[v_A]$. So
$$\sh_A=\ov{ \Slt (\C[t])\cdot [v_A]}\hk\Pro (M^A).$$
These are complex projective algebraic varieties.
They can be described as a subvarieties of the product of the affine
$\slt$ grassmanians. In fact, this is a consequence from the results
from \cite{mi2}, where the fusion products were embedded into the
tensor product of the integrable irreducible level $1$ $\slth$ modules.
Let us give some details.

For $A=(a_1\le\ld\le a_n)\in\N^n$ define $D=(d_1,\ld,d_{a_n})$ with
$d_i=\#\{j:\ a_j=i\}$. In \cite{mi2} for
any such $D$ the
integrable $\slth$ module $L^D$ was constructed. It has two descriptions:
as an inductive limit of some special fusion products
and as a submodule of the tensor product of the irreducible level one $\slth$
modules. The first construction is the following one.

Let $a_i=a_{i+1}$. In \cite{mi2} the
following exact sequence of $\sltc$ modules was constructed:
\begin{equation}
\label{sim}
0\to M^{(a_1,\ld,a_{i-1},a_{i+2},\ld,a_n)}\to M^A\to
M^{(a_1,\ld,a_{i-1},a_i-1,a_{i+1}+1,a_{i+1},\ld,a_n)}\to 0
\end{equation}
(note that the corresponding relation for the $q$-characters of the
modules from $(\ref{sim})$ can be found in \cite{SW1, SW2}).
Introduce a notation
$$A^{(i)}=(a_1,\ld,a_n,\underbrace{a_n,\ld,a_n}_{2i}).$$
From $(\ref{sim})$ we obtain an embedding
$M^{A^{(i)}}\hk M^{A^{(i+1)}}$. This allows us to construct an inductive limit
of the modules $M^{A^{(i)}}$, which is denoted as $L^D$:
\begin{equation}
\label{indlim}
L^D=M^{A}\hk M^{A^{(1)}} \hk M^{A^{(2)}}\hk\ld.
\end{equation}
By its definition, $L^D$ is $\sltc$
module. It was
proved in \cite{mi2} that in fact it carries a structure of the level
$(a_n-1)$ $\slth$ module.

These modules can be also described as a submodules of the tensor product
of the level one $\slth$ modules.
Let $L_{0,1}$ and $L_{1,1}$ be two level one irreducible integrable
$\slth$ modules ($L_{0,1}$ is a vacuum module). Denote by $v(2k)$ the set
of the extremal vectors in $L_{0,1}$ and by $v(2k+1)$ the set of the
extremal vectors in $L_{1,1}$. For example, $v(0)\in L_{0,1}$ is a vacuum
vector, $v(1)\in L_{1,1}$ is a highest vector and
$$e_{i-1}v(i)=v(i-2);\qquad e_k v(i)=0 \text{ if } k>i-1.$$
We have
$$L^D=U(\slth)\cdot
\left(v(n-d_1)\T v(n-d_1-d_2)\T\ld \T v(n-d_1-\ld-d_{a_n-1})\right).$$
This module is a submodule of the tensor product
$L_{p_1,1}\T\ld\T L_{p_{a_n-1},1}$, where $p_i=0$, if $n-d_1-\ld-d_{i}$ is
even and $p_i=1$ otherwise.

Now we can embed the Schubert varieties $\sh_A$ into the product of the
affine $\slt$ grassmanians.
Let $Gr_0$ and $Gr_1$ be even and odd
$\slt$ affine grassmanians, i.e. a subvarieties of
$\Pro (L_{0,1})$ and $\Pro (L_{1,1})$ correspondingly:
$$Gr_0=\Slth\cdot [v(0)],\qquad Gr_1=\Slth \cdot [v(1)].$$
Let $\Ob_A$ be an $\Slth$-orbit in the product
$Gr_{p_1}\times \ld \times Gr_{p_{a_n-1}}$:
$$\Ob_A=\Slth\cdot
\left([v(n-d_1)]\times [v(n-d_1-d_2)]\times\ld \times
[v(n-d_1-\ld-d_{a_n-1}]\right).$$
Then we define the generalized affine $\slt$ grassmanian as the closure
of $\Ob_A$ in $Gr_{p_1}\times \ld \times Gr_{p_{a_n-1}}$:
$Gr_A=\ov{\Ob_A}$.
Surely, $Gr_A$ is a subvariety of the projective
space $\Pro (L^D)$ (again, $Gr_A$ is a closure of the $\Slth$-orbit of the
point).
Let us give some examples:
\begin{enumerate}
\item $Gr_{(2,3)}\hk Gr_0\times Gr_1$ is the flag manifold for the
group $\Slth$.
\item If $A=(\underbrace{2,\ld,2}_{2n})$, then $Gr_A=Gr_0$.
\item If $A=(\underbrace{2,\ld,2}_{2n+1})$, then $Gr_A=Gr_1$.
\end{enumerate}

Because of the definition $(\ref{indlim})$ we have the
approximation of $Gr_A$ by the finite-dimensional
Schubert varieties.
Namely, the generalized grassmanian is an inductive limit
$$Gr_A=\sh_A\hk \sh_{A^{(1)}}\hk \sh_{A^{(2)}}\hk\ld.$$
(For the another approach to the study of the
infinite-dimensional varieties of the above type see \cite{sto1}).

The structure of the Schubert varieties $\sh_A$ depends on the number
of the different elements among $a_i$. For example, if all the numbers $a_i$
are different, then $\sh_A$ is smooth. Otherwise, it has singularities, and
the structure of these singularities depends on the number of such $i$ that
$a_i<a_{i+1}$ (recall that we put $a_i\le a_{i+1}$). In this paper we
concentrate on the study of the smooth case (though some theorems will be
proved in the whole generality).

We start with a proof that for all
$A$ with $a_i< a_{i+1}$ the corresponding Schubert varieties are isomorphic.
In order
to prove that we construct a bundle $\sh_A\to\Pro^1$ with a fiber
$\sh_{(a_1,\ld,a_{n-1})}$ and show that the transition functions
of this bundle do not depend on $A$. Let us give some details.
Let $u_A$ be
the highest weight vector in $M^A$ (so $u_A=w\cdot v_A$, where $w$ is the
Weyl element from $\Slt\hk \Slt (\C[t]/t^n)$). There are two so-called
big cells in $\sh_A$:
$$U_x=\exp\left(\sum_{i=0}^{n-1} e_i x_i\right)\cdot [v_A],\quad
  U_y=\exp\left(\sum_{j=0}^{n-1} f_j y_j\right)\cdot [u_A],\qquad
  x_i,y_j\in\C.$$
Both of these cells are isomorphic to $\C^n$ and the closure of each is
$\sh_A$. We also have that $U_x\cup U_y=\G \cdot [v_A]$,
where $\G=\Slt (\C[t]/t^n)$. Moreover,
$$\G\cdot [v_A]=
U_x\bigcup \exp\left(\sum_{i=1}^{n-1} f_iy_i\right)\cdot [u_A].$$
We prove that
$\exp\left(\sum_{i=0}^{n-1} e_i x_i\right)\cdot [v_A]=
\exp\left(\sum_{j=0}^{n-1} f_j y_j\right)\cdot [u_A]$ if
\begin{equation}
\label{trfun}
(x_0+tx_1+\ld+t^{n-1}x_{n-1})(y_0+ty_1+\ld+t^{n-1}y_{n-1})=1
\end{equation}
in the ring $\C[t]/t^n$.
Now, one can define a map $U_x\to\Slt\cdot[v_A]\simeq\Pro^1$ by the formula
$\exp\left(\sum_{i=0}^{n-1} e_i x_i\right)\cdot [v_A]\to
\exp(e_0x_0)\cdot [v_A]$. We prove that this map can be extended to the
bundle $\shn\to \Pro^1$ with the fiber $\sh^{(n-1)}$ (this fiber comes
from the following statement -- see \cite{mi1}:
$\C[e_1,\ld,e_{n-1}]\cdot [v_A]\simeq M^{(a_1,\ld,a_{n-1})}$).
So from the formula $(\ref{trfun})$ we obtain
that $\sh_A$ is independent on $A$. We denote this variety by $\shn$.

The latter can be described as a generalized flag manifold ${\bf F}$
(see \cite{segal}).
Consider the variety of the sequences of the spaces
$$\C[t]\oplus \C[t]=
V_0\hookleftarrow V_1 \hookleftarrow \ld \hookleftarrow V_n$$
with a properties
$tV_i\hk V_{i+1}$ and $\dim (V_{i-1}/V_i)=1$.
Using the obvious projection ${\bf F}\to \Pro^1\simeq \Pro(V_0/tV_0)$
and an action of
$\Slt(\C[t])$ on ${\bf F}$ we prove that ${\bf F}$ is isomorphic to
$\shn$. (Note that for any $A$ with possibly coinciding elements the variety
$\sh_A$ can be also described via the construction of this type,
but we do not study it in the present paper).

We have already mentioned that there exists a bundle
$\shn\to \Pro^1$ with a fiber $\sh^{(n-1)}$.
In fact,
we have even more interesting picture: we prove that there are in-between
bundles $\shn\to \sh^{(n-1)}\to\ld\to \sh^{(2)}\to \sh^{(1)}$ (note that
$\sh^{(1)}=\Pro^1$), such that any map $\sh^{(m+k)}\to \sh^{(m)}$ is a
bundle with a fiber $\sh^{(k)}$. We denote such maps by $\pi_{m+k,m}$.

The next point of the study is the Lie algebra of vector fields on $\shn$.
We prove that it practically coincides with a
"base" Lie algebra $\slt\T (\C[t]/t^n)$.
There are some additional operators, acting on $M^A$: these are an
annihilation operators $L_0,\ld,L_{n-2}$ from the Virasoro algebra
(for $x\in\slt$ we have $[L_i,x\T t^j]=-jx\T t^{i+j}$).
As we have already mentioned, $M^A$ is embedded into the integrable
$\slth$ module $L^D$. Thus, we have an action of the Virasoro algebra
on $L^D$ and hence on $M^A$.
Note that the existence of the vector fields from $\slt\T (\C[t]/t^n)$
is obvious from the definition of $\shn$. It is easy to show that Virasoro
operators also
generates the vector fields on $\shn$. We prove that
$e_i, h_i, f_i, \ i=0,\ld,n-1$ (here $x_i=x\T t^i$) and $L_i, i=0,\ld,n-2$
form a basis of the Lie algebra $\Vect(\shn)$. For the proof, we compute the
transition matrix for the bundle of the vector fields, which are tangent
to the fibers of $\pi_{n,1}$. This allows us to show that
$\dim \Vect (\shn)=4n-1$.

Note that the first definition of the varieties $\sh_A$ was given via
the abelian Lie algebra with no $\slt$ mentioned.
But the current algebra $\slt\T \C[t]$ naturally appears as a main part
of the Lie algebra of the vector fields on the Schubert variety.

In \cite{mi2} some natural submodules of the fusion product were mentioned.
In this paper we study the structure of these submodules.
They are important because of
their  connection with some
$\Slt (\C[t])$-invariant subvarieties of
$\shn$. Let us give some details.

Using the properties of the fusion products
(see \cite{mi1,mi2}) one can show that for any $i=1,\ld,n-1$ there exists
a submodule $S_{i,i+1}(A)\hk M^A$, such that the following sequence
of $\sltc$ modules is exact:
\begin{equation}
\label{submodules}
0\to S_{i,i+1}(A)\to M^A\to M^{(a_1,\ld,a_{i-1},a_i-1,a_{i+1}+1,\ld,a_n)}
\to 0.
\end{equation}
(Note that we have already used the special case of this sequence in
$(\ref{sim})$).
For example, if $i=1$ then $S_{i,i+1}(A)\simeq M^{(a_2-a_1+1,a_3,\ld,a_n)}$
(the corresponding recurrent formula for the $q$-characters can be found in
the works of Shilling and Warnaar). In the case $i>1$ the structure of
$S_{i,i+1}(A)$ is more complicated.
We prove, that if
$i>1, a_{i-1}>1$ and $a_i\ne a_{i+1}$, then there exists an exact sequence
\begin{multline}
0\to M^{(a_1,\ld,a_{i-2},a_{i-1}-a_i+a_{i+1},a_{i+2},\ld,a_n)}\to
S_{i,i+1}(A)\to\\ \to
S_{i,i+1}(a_1,\ld,a_{i-2},a_{i-1}-1,a_i+1,a_{i+1},\ld,a_n)\to 0.
\end{multline}

Another description of the submodules $S_{i,i+1}(A)$ is the following.
Denote
$$A'=(a_1,\ld,a_{i-1},\underbrace{a_i,\ld,a_i}_{n-i-1});\qquad
A''=(a_{i+1}-a_i+1, a_{i+2}-a_i+1,\ld, a_n-a_i+1).$$
Recall that $v_{A'}$ and $v_{A''}$ are the cyclic vectors of $M^{A'}$ and
$M^{A''}$.
We prove that there is an embedding of $\slt\T (\C[t]/t^n)$ modules
$S_{i,i+1}(A)\hk M^{A'}\T M^{A''}$ and $S_{i,i+1}(A)$ is generated from the
vector
$v_{A'}\T v_{A''}$ by the action of the polynomials in the variables
$e_j,\ j=0,\ld, n-3$ and $\Id\T e_{n-i-1}$, where $\Id$ is acting on
$M^{A'}$ and $e_{n-i-1}$ on $M^{A''}$.
Thus we obtain an additional operators, acting on $S_{i,i+1}(A)$ (surely, we
have operators from $\slt\T (\C[t]/t^n)$). The Lie algebra, acting on
$S_{i,i+1}(A)$ can be described as follows.
Let $\Lie_{i,n},$ $0< i\le n$ be a graded Lie algebra,
$\Lie_{i,n}=\slt\T B_{i,n},$ where $B_{i,n}$ is a commutative graded
associative algebra with generators $t$ of degree $1$ and $u$ of degree
$n-i$ and relations $t^n=0$, $tu=0$, $u^2=0$.
For $i=0$ let
$\Lie_{0,n}$ be $\slt\T (\C[t]/t^{n+1}).$
We prove that
$S_{i,i+1}(A)$ is cyclic $\Lie_{i-1,n-2}$ module with a cyclic vector being
the tensor product $v_{A'}\T v_{A''}=v_{i,i+1}(A)$.
This allows to make a conjecture (we do not discuss it in the paper) that
$S_{i,i+1}(A)$ is a fusion product
itself, but with respect to the algebra $\Lie_{i-1,n-2}$. It means that
for some filtration on the Lie algebra $\bigoplus_{i=1}^{n-1} \slt$, with an
adjoint graded algebra being $\Lie_{i-1,n-2}$, the induced filtration on
the tensor product of $\slt$ modules
$$\C^{a_1}\T\ld \C^{a_{i-1}}\T\C^{a_{i+1}-a_i+1}\T \C^{a_{i+2}}\T \ld \T
\C^{a_n}$$
gives $S_{i,i+1}(A)$ as an adjoint graded module.

Now let us return to the varieties $\shn$. Recall that there is a "big"
dense orbit $\G \cdot [v_A]\hk \shn$.
We want
to understand the structure of the complement
$\shn\setminus \G \cdot [v_A]$.
We prove that this complement is the union of $n-1$ irreducible varieties of
complex dimension $n-1$, which are
denoted by $\shin, i=1,\ld,n-1$. They can be described as follows:
$\shin$ is the closure of the
orbit $\Lg_{i-1,n-2}\cdot [v_{i,i+1}(A)]$ in the projective space
$\Pro(S_{i,i+1}(A))$, where $\Lg_{i-1,n-2}$ is the Lie group of the Lie
algebra $\Lie_{i-1,n-2}$. In fact, the following is true:
$\shin=\shn\cap \Pro(S_{i,i+1}(A))$.

We can give another description of $\shin$ as a
subvariety of the product of two Schubert varieties.
Namely, we prove that
$$\sh_i^{(n-1)}=\{(x,y)\in \sh^{(n-2)}\times \sh^{(n-i)}:\
\pi_{n-2,n-i-1}(x)=\pi_{n-i,n-i-1}(y)\}.$$
For example, $\sh^{(n-1)}_1\simeq \sh^{(n-1)}$ and
$\sh^{(n-1)}_{n-1}\simeq \sh^{(n-2)}\times \Pro^1$.
Note also that the classes of $\shin$, $i<n$ and the class of the fiber
of the map $\pi_{n,1}$
are the generators of the group $H_{2n-2}(\shn,\Z)$.

One of the goals of this paper is
to realize the fusion products as the dual spaces
of the sections of some line bundles on $\shn$.
Let
$$C_i=\ov{\{\exp(t e_i)\cdot [v_A], t\in\C\}}, i=0,\ld,n-1.$$
(In particular,
$C_0=\Slt\cdot [v_A]$). All these curves are $\Pro^1$ and their fundamental
classes are the generators of the group $H_2(\shn,\Z)$.
Note also that the classes of $C_{n-i}$, $0< i\le n$
form the dual basis with respect to the basis of the group
$H_{2n-2}(\shn,\Z)$, mentioned above.
One can show that any line bundle on $\shn$ is uniquely
determined by its restriction to the lines $C_i$.
Let $\E$ be a line bundle and $\E|_{C_i}=\Ob(b_i)$.
We denote such bundle as $\Ob(b_{n-1}, b_{n-2}-b_{n-1},\ld,b_0-b_1)$.
\begin{rem*}
Recall that we have already mentioned that there exists a variety
$\Sh_A$, which is fibered over $\C^n$ with a fiber $\sh_A(Z)$. If
$a_i\ne a_j$, then $\Sh_A$ doesn't depend on $A$ and we denote it by
$\Sh_n$. Note that any line bundle $\Ob(c_1,\ld,c_n)$ on $\sh_n$ is
a restriction of some line bundle on $\Sh_n$ to the "special" fiber.
The
restriction of the latter bundle to the general fiber -- which is simply
$\Pro^1\times\ld\times \Pro^1$ -- is
$\Ob(c_1)\boxtimes\ld\boxtimes \Ob(c_n)$.
\end{rem*}

Recall that for all $A$ with pairwise distinct elements we have an embedding
$\imath_A: \shn\hk \Pro(M^A)$. We prove that for exceptional $A$ (with
$a_i=a_j$ for some $i\ne j$) there exists a surjective birational map
$m_A:\shn\to\sh_A$
and hence a map $\widehat{\imath_A}: \shn\to \Pro(M^A)$, which is
a composition of $m_A$ and the embedding of $\sh_A$ to $\Pro(M^A)$ (surely,
$\widehat{\imath_A}$ is not an embedding). Note that in fact the map
$m_A$ is a resolution of the singularities of the singular variety $\sh_A$.
A natural way of constructing the line bundles on $\shn$ is taking the
inverse image of the bundles on $\Pro(M^A)$.
It is easy to show that
$$\imath_A^* \Ob(1)=\Ob(a_1-1,\ld,a_n-1)$$
and for the exceptional $A$
$$\widehat{\imath_A}^* \Ob(1)=\Ob(a_1-1,\ld,a_n-1).$$
The main theorem, proved in the last section, states that if
$1\le a_1\le a_2\le \ld \le a_n$,
then there is an isomorphism of
$\slt\T \C[t]/t^n$ modules:
\begin{equation}
\label{intsec}
H^0(\shn, \Ob(a_1-1,\ld,a_n-1))^*\simeq M^A.
\end{equation}
Because of the existence of the maps $\imath_A$ and $\widehat{\imath_A}$
it is enough to show that
$$\dim H^0(\shn, \Ob(\an))=\prod_{i=1}^n (a_i+1).$$
We prove a little bit more general statement: let $0\le a_1\le\ld\le a_n$.
Then
the dimension of the zeroth cohomologies of the bundle $\Ob(\an)$ is equal
to $\prod_{i=1}^n (a_i+1)$  and all the higher
cohomologies vanish. For the proof we restrict our bundles to the divisors
$\shin$ and write the corresponding exact sequences.
Note also that it is interesting to compute all cohomologies for all
line bundles on $\shn$ (recall that we have an answer only for
$\Ob(\an)$ with nonnegative and nondecreasing $a_i$). We have a conjecture
that in the general case the structure of the cohomologies as an
$\slt\T\C[t]$ modules can be given in terms of the fusion products.

Now let us clarify the connection between the exact sequence
$(\ref{submodules})$
and the subvarieties $\shin$. Let $0\le a_1\le \ld \le a_n$. Let $J_i$ be
the subspace of the sections of the bundle $\Ob(\an)$:
$$J_i=\{ s\in H^0(\shn,\Ob(\an)):\ s|_{\shin}=0\}.$$
Then as $\sltc$ module the dual space $J_i^*$ is isomorphic to the fusion
product
$M^{(\ld,a_i-1,a_{i+1}+1,\ld)}$ and we have an exact sequence
(which is dual to $(\ref{submodules})$:
$$0\to J_i\hk H^0(\shn,\Ob(\an))\sur S_{i,i+1}(A)^*\to 0.$$
This gives us the geometric description of the modules $S_{i,i+1}(A)$:
if $\jmath$ is an embedding $\shin\hk\shn$, then
$$S_{i,i+1}(A)^*\simeq H^0(\shin,\jmath_*\jmath^*\Ob(\an)).$$

In the end, let us return to the infinite-dimensional varieties $Gr_A$. As we
have already mentioned, these generalized grassmanians are the inductive limit
of the finite-dimensional Schubert varieties $\sh_{A^{(i)}}$, where
$$A^{(i)}=(a_1,\ld,a_n,\underbrace{a_n,\ld,a_n}_{2i}).$$
Thus, for the study of the bundles on $Gr_A$ we need to study the singular
varieties $\sh_{A^{(i)}}$. Although we are not doing it in this paper, let us
point the main moments. There exists a line bundle $\Ob$ on $Gr_A$ such that
the
dual space $H^0(\Ob)^*$ is isomorphic to the module $L^D$ (there
$D=(d_1,\ld,d_{a_n})$, $d_i=\#\{j:\ a_j=i\}$) (this is a consequence of the
corresponding statement for the finite-dimensional Schubert varieties).
Recall that in \cite{mi2} the module $L^D$ was decomposed into the
irreducible components and the multiplicity $c_i$
of each irreducible module $L_{i,a_n-1}$ in the decomposition was given in
terms of the
Verlinde algebra. Moreover, each $c_i$ inherit a $q$-grading from $L^D$.
One can show that the corresponding polynomial in $q$ for $c_0$ is a
restricted Kostka polynomial for $\slt$ (see \cite{kost1,kost2}).

The paper is organized in the following way:

In the first section we study the geometric properties of $\shn$. We prove,
that there is a bundle $\shn\to\Pro^1$ with a fiber $\sh^{(n-1)}$
(theorem $(\ref{bunth})$) and a more general fact that for any $m>k$ there
exists a bundle $\sh^{(m)}\to\sh^{(k)}$ with a fiber $\sh^{(m-k)}$
(corollary $(\ref{inbet})$). We also compute the Lie algebra of the vector
fields on $\shn$. The answer is given in the theorem $(\ref{vectf})$.

The second section is devoted to the study of the structure of the fusion
product. Namely, the submodules $S_{i,i+1}(A)$ are studied. Three
descriptions of the latter are given. First we construct a filtration
with an adjoint quotients being some fusion products (proposition
$(\ref{filt})$). The second and the most important (for the geometry of the
varieties $\shin$) description gives the embedding of the module
$S_{i,i+1}(A)$ to the tensor product of some special fusion products
(proposition $(\ref{mprop})$). And in the lemmas $(\ref{ind1})$,
$(\ref{ind2})$,
we give an inductive description of $S_{i,i+1}(A)$ via
$S_{i,i+1}(a_1,\ld,a_{n-1})$.

In the last section we continue the study of the geometric properties of the
Schubert varieties. First in the proposition $(\ref{coordring})$ we prove,
that $\shn$ is a projective algebraic
variety: we compute its coordinate ring
(for the similar "functional" construction see
\cite{Gordon}).
Then we
give the description of the varieties $\shin$ in terms of the Schubert
varieties (proposition $(\ref{descr})$). The paper finishes with the theorem
$(\ref{maintheorem})$, which computes the cohomologies of the bundles
$\Ob(\an)$ $(0\le a_1\le\ld\le a_n)$. As a corollary we obtain the
realization
of the fusion products in the dual space of sections of the line bundles
(corollary $(\ref{sections})$). \\
{\bf Acknowledgements.}
The first named author was partially supported by the grants
SS 2044.2003.2, INTAS 00-55.
The second named author was partially supported by
the RFBR grant 03-01-00167.

\section {Geometric structure of the Schubert varieties}
\subsection {Definition and first properties}
Let $e,h,f$ be a standard $\slt$ basis.
In \cite{fusion} the set of
$\sltc$ modules $M^A$, $A=(\an)\in\N^n$, called fusion products was defined.
These modules are also denoted as $\Ca$.
We briefly recall the definition and main properties (see \cite{fusion}).

$M^A$ is an adjoint graded module with respect to the below filtration
of the tensor product of
the evaluation $\sltc$ modules
$$\C^{a_1}(z_1)\T\ld\T\C^{a_n}(z_n).$$
Here $z_i$ are pairwise distinct complex numbers and filtration of the
tensor product is the following:
$$F_s=\bra e_{i_1}\ld e_{i_k}\ov{v_A},\ \sum i_j\le s\ket$$
($x_i=x\T t^i, x\in\slt$,\  $\ov{v_A}=v_1\T\ld\T v_n$ and $v_i$ is a
lowest vector of $\C^{a_i}(z_i)$).
Thus, $M^A$ is an $\Cn$ cyclic module with a cyclic
vector $v_A$ -- the image of $\ov{v_A}$:
$$M^A=\Cn/I_A$$
and  $I_A$ is some ideal in the polynomial ring $\Cn$.
(Surely, $M^A$ is also $\slt\T\C[t]$ module).
In \cite{mi1} it
is shown that $I_A$ is generated by the elements
\begin{equation}
\label{rel}
\sum_{i_1+\ld +i_k=s} e_{n-1-i_1}e_{n-1-i_2}\ld e_{n-1-i_k},\ k\ge 1,\
s<\sum_{j=1}^n (k+1-a_j)_+
\end{equation}
(here $a_+=\max (a,0)$).

Denote also by $u_A\in M^A$ the image of the product of the highest weight
vectors of
$\C^{a_i}(z_i)$. Then $M^A$ is also a cyclic $\C[f_0,\ld,f_{n-1}]$ module
with a cyclic vector $u_A$ and relations of the form $(\ref{rel})$.

Consider a Lie group $\G=\Sltc$ of a Lie algebra $\slt\T (\C[t]/t^n)$.
$\G$ acts on $M^A$ and thus on its projectivization $\Pro (M^A)$.
For $v\in M^A$ we denote by $[v]\in \Pro (M^A)$  the line $\C v$.
\begin{opr}
Let $A=(a_1,\ld,a_n)\in\N^n$. The Schubert variety $\sh_A$ is a closure
of an orbit of the point $[v_A]$ in $\Pro (M^A)$:
$$\sh_A=\overline{\G\cdot [v_A]}.$$
\end{opr}

Denote by $\U_+$ and $\U_-$ the subgroups of $\G$, generated by
$\exp(p_ie_i)$ and $\exp(q_if_i)$ respectively ($p_i,q_i\in\C, i=0,\ld,n-1$).
Thus, both $\U_+$ and $\U_-$ are isomorphic to $\C^n$,
\begin{equation}
\label{unip}
\U_+=
\begin{pmatrix}
1 & p(t) \\ 0 & 1
\end{pmatrix}\text{,}\ \ \
\U_-=
\begin{pmatrix}
1 & 0 \\ q(t) & 1
\end{pmatrix}{,}\ \ \ \
p(t),q(t)\in\C[t]/t^n.
\end{equation}

Now we will formulate some lemmas, which can be easily checked.
\begin{lem}
We have an isomorphism:
$$\U_+\cdot [v_A]\simeq\C^n;\ \ \
\exp\left(\sum_{i=0}^{n-1} p_ie_i\right)\cdot [v_A]
\mapsto (p_0,\ld,p_{n-1}).$$
We call this set a big cell.
(We have also another big cell: $\U_-\cdot [u_A]\simeq\C^n).$
\end{lem}
Denote
$
w=
\begin{pmatrix}
0  & 1\\
-1 & 0
\end{pmatrix}
\in\Slt\subset\G.
$
Note that $w\cdot [v_A]=[u_A]$.
Let $\U^{(i)}_+$ and $\U^{(i)}_-$ be a subgroups of $\U_+$ and $\U_-$:
$$ \U^{(i)}_+=
\begin{pmatrix}
1 & t^ip(t) \\ 0 & 1
\end{pmatrix},\ \ \
\U^{(i)}_-=
\begin{pmatrix}
1 & 0 \\ t^iq(t) & 1
\end{pmatrix}
.
$$
\begin{lem}
\label{strat}
\begin{enumerate}
\item $\overline {\U_+\cdot [v_A]}=\sh_A=\ov {\U_-\cdot [u_A]}.$
\item $\G\cdot [v_A]=\U_+\cdot [v_A]\bigsqcup \U^{(1)}_-\cdot [u_A].$
\end{enumerate}
\end{lem}

\begin{lem}
\label{bunlem}
$1$.\ $\Slt\cdot [v_A]\simeq \Pro^1$.\\
$2$.\ There is an $\Slt$-equivariant bundle
$\widetilde{\pi}: \G\cdot [v_A]\to \Slt\cdot [v_A]$ with a fiber $\C^{n-1}$:
$$\exp\left(\sum_{i=0}^{n-1} p_ie_i\right) [v_A] \mapsto
\exp(e_0p_0) [v_A];\
\exp\left(\sum_{i=1}^{n-1} q_if_i\right) [u_A] \mapsto [u_A]$$
$3$.\ Let $U_0$ and $U_1$ be two charts of $\Pro^1\simeq\Slt\cdot [v_A]$:
$$U_0=\Pro^1\setminus [u_A],\ U_1=\Pro^1\setminus [v_A].$$
Then
$$\widetilde{\pi}^{-1} U_0=\U_+\cdot [v_A] \simeq \{p(t)\in\C[t]/t^n\},\
  \widetilde{\pi}^{-1} U_1=\U_-\cdot [u_A] \simeq \{q(t)\in\C[t]/t^n\},$$
and two polynomials $p(t)$ and $q(t)$ correspond to the same
point of $\G\cdot [v_A]$ if $p(t)q(t)=1$ in $\C[t]/t^n$.
\end{lem}
\begin{proof}
The last statement follows from the equality in the group $\G$:
$$
\begin{pmatrix}
1 & p(t)\\
0 & 1
\end{pmatrix}
\begin{pmatrix}
p(t) & 0\\
-1 & p(t)^{-1}
\end{pmatrix}
=
\begin{pmatrix}
1 & 0\\
p(t)^{-1} & 1
\end{pmatrix}
\begin{pmatrix}
0 & 1\\
-1 & 0
\end{pmatrix}
$$
\end{proof}

Let $A=(a_1,\ld,a_n)$, $a_1\le\ld\le a_n$.
The natural question is: is it possible to extend our bundle to
the closure of the orbit $\G\cdot [v_A]$? The answer is positive
if $a_{n-1}<a_n$ and negative otherwise. The impossibility of the extension
means that in the latter case the closures of the fibers in the different
points intersect. For the proof we need a lemma, which was proved in
\cite{mi1}.
\begin{lem}
\label{demazure}
Let $\varphi$ be an isomorphism of the algebras:
$$\varphi: \C[e_1,\ld,e_{n-1}]\to \C[e_0,\ld,e_{n-2}],\qquad
e_i\mapsto e_{i-1}.$$
Let $I^{(1)}_A$ be an ideal in $\C[e_1,\ld,e_{n-1}]$ such that
$$\C[e_1,\ld,e_{n-1}]\cdot v_A\simeq \C[e_1,\ld,e_{n-1}]/I^{(1)}_A$$
(i.e. $$I^{(1)}_A=\{f(e_1,\ld,e_{n-1}):\ f(e_1,\ld,e_{n-1})v_A=0\}).$$
Then $\varphi (I^{(1)}_A)=I_{(a_1,\ld,a_{n-2})}$.
In other words
$$\C[e_1,\ld,e_{n-1}]\cdot v_A\simeq
\C[e_0,\ld,e_{n-2}]\cdot v_{(a_1,\ld,a_{n-1})}.$$
We also have an isomorphism of $\C[f_0,\ld,f_{n-1}]$ modules
$$M^A/\C[e_1,\ld,e_{n-1}]\cdot v_A\simeq M^{(a_1,\ld,a_{n-1},a_n-1)}.$$
\end{lem}

\begin{theorem}
\label{bunth}
Let $A=(\an)\in\N^n,\ a_1\le \ld\le a_{n-1}<a_n$.
Then the bundle $\widetilde{\pi}: \G\cdot [v_A]\to\Pro^1$ can be extended
to the
bundle $\pi: \sh_A\to \Pro^1$ with a fiber $\sh_{(a_1,\ld,a_{n-1})}$.
For any $x\in\Pro^1$ we have
$$\pi^{-1}(x)=\ov{\widetilde{\pi}^{-1}(x)}.$$
\end{theorem}
\begin{proof}
Because of the lemma $(\ref{demazure})$
the closure of $\widetilde{\pi}^{-1}([v_A])$ in $\sh_A$ is
isomorphic to $\sh_{(a_1,\ld,a_{n-1})}$. So, the only thing to be proved
is a fact that the closures of the fibers of $\widetilde{\pi}$ in the
different points are isomorphic and do not intersect.

Note that for any $z\in\C$
$$\ov{\widetilde{\pi}^{-1}(\exp (ze_0)[v_A])}\subset
\Pro (\C[e_1,\ld,e_{n-1}]\cdot (\exp(ze_0)v_A)).$$
Hence, because of the action of the group $\Slt$ on $\shn$ it is enough to
prove that
\begin{equation}
\label{eq}
\C[e_1,\ld,e_{n-1}]\cdot (\exp(ze_0)v_A)\simeq \C[e_1,\ld,e_{n-1}]\cdot v_A
\end{equation}
and
\begin{equation}
\label{int}
\C[e_1,\ld,e_{n-1}]\cdot (\exp(ze_0)v_A) \bigcap
\C[e_1,\ld,e_{n-1}]\cdot v_A=0.
\end{equation}

The statement $(\ref{eq})$ is obvious, because the operator $\exp(ze_0)$ is
invertable and thus provides an isomorphism between the left and right
hand sides. Note that both left and right hand sides are isomorphic to
$M^{(a_1,\ld,a_{n-1})}$.

Let us prove the statement $(\ref{int})$.
First note that if
$$\C[e_1,\ld,e_{n-1}]\cdot v_A\bigcap
\C[e_1,\ld,e_{n-1}]\cdot (\exp(ze_0)v_A)\ne 0$$
then $\C[e_1,\ld,e_{n-1}]\cdot v_A$ and $\C[e_1,\ld,e_{n-1}]\cdot (e_0v_A)$
also have a nontrivial intersection. In fact, let $p(e_1,\ld,e_{n-1})$ be
a homogeneous polynomial and
$$0\ne p(e_1,\ld,e_{n-1})(\exp (ze_0)v_A)\in \C[e_1,\ld,e_{n-1}]\cdot v_A.$$
Then either $p(e_1,\ld,e_{n-1})(e_0 v_A)$ is a nontrivial element of
$\C[e_1,\ld,e_{n-1}]\cdot v_A$ or
$$p(e_1,\ld,e_{n-1})(e_0 v_A)=0\qquad \text{ and }\qquad
p(e_1,\ld,e_{n-1}) v_A\ne 0.$$
Let us show that the last variant is impossible, i.e. that
$$\C[e_1,\ld,e_{n-1}]\cdot v_A\simeq \C[e_1,\ld,e_{n-1}]\cdot (e_0v_A).$$
Because of the lemma $(\ref{demazure})$ we have an isomorphism
$$M^A/\C[e_1,\ld,e_{n-1}]\cdot v_A\simeq M^{(a_1,\ld,a_{n-1},a_n-1)}$$
with the cyclic vector of the right hand side being the image of $e_0v_A$.
But in $M^{(a_1,\ld,a_{n-1},a_n-1)}$ we have an isomorphism
$$\C[e_1,\ld,e_{n-1}]\cdot v_{(a_1,\ld,a_{n-1},a_n-1)}\simeq
M^{(a_1,\ld,a_{n-1})}$$
(because $a_{n-1}<a_n$ !). Thus in the quotient
$M^A/\C[e_1,\ld,e_{n-1}]\cdot v_A$ the image of
$\C[e_1,\ld,e_{n-1}]\cdot (e_0v_A)$  is isomorphic to $M^{(a_1,\ld,a_{n-1})}$
and so to $\C[e_1,\ld,e_{n-1}]\cdot v_A$.

Now the only thing left is to prove that
$$\C[e_1,\ld,e_{n-1}]\cdot v_A\bigcap \C[e_1,\ld,e_{n-1}]\cdot (e_0v_A)=0.$$
But that is a consequence of the already used fact that in
$M^A/\C[e_1,\ld,e_{n-1}]\cdot v_A$ the image of
$\C[e_1,\ld,e_{n-1}]\cdot (e_0v_A)$ is isomorphic to $M^{(a_1,\ld,a_{n-1})}$.
Theorem is proved.
\end{proof}

\begin{cor}
Let $A=(a_1,\ld,a_n),\ a_1<a_2<\ld <a_n$. Then $\sh_A$ is a smooth
$n$-dimensional complex manifold.
\end{cor}
\begin{proof}
We know that $\sh_A$ is fibered over $\Pro^1$ with a fiber
$\sh_{(a_1,\ld,a_{n-1})}$, which is smooth by the induction.
\end{proof}

\begin{cor}
Let $A=(a_1<\ld <a_n), B=(b_1<\ld <b_n)$. Then $\sh_A$ is isomorphic to
$\sh_B$.
\end{cor}
\begin{proof}
We use the induction on $n$.
$\sh_A$ and $\sh_B$ are fibered over $\Pro^1$ with a fibers, which are
isomorphic by the induction assumption. Moreover, because of the lemma
$(\ref{bunlem})$ we know that $\G\cdot [v_A]\simeq\G\cdot [v_B]$,
because the
transition functions of the bundle $\G\cdot [v_A]\to\Pro^1$ do not depend on
$A$. But $\ov{\G\cdot [v_A]}=\sh_A$ and $\ov{\G\cdot [v_B]}=\sh_B$.
Thus we obtain $\sh_A\simeq \sh_B$.
\end{proof}

Introduce a notation
$$\sh^{(n)}=\sh_A,\ \ \ A=(a_1<a_2,\ld <a_n).$$
Note that $\sh^{(1)}\simeq \Pro^1$. Thus, the bundle $\shn\to \Pro^1$
from the theorem $(\ref{bunth})$ can be regarded as a bundle
$\shn\to \sh^{(1)}$ with a fiber $\sh^{(n-1)}$.
We will prove that there are in-between bundles:
\begin{equation}
\label{map}
\shn\to  \sh^{(n-1)}\to\ld\to\sh^{(2)}\to\sh^{(1)}
\end{equation}
with a fibers $\Pro^1\simeq\sh^{(1)}$.
Moreover, any composition of the maps from $(\ref{map})$ is a bundle
$\sh^{(m+k)}\to \sh^{(m)}$ with a fiber $\sh^{(k)}$.
For the proof we need the following proposition:
\begin{prop}
\label{tg}
Let $A=(a_1\le\ld\le a_n), B=(b_1\le\ld\le b_m),\ n\ge m$. Denote
$$C=(a_1,\ld,a_{n-m},a_{n-m+1}+b_1-1,a_{n-m+2}+b_2-1,\ld,a_n+b_m-1).$$
Consider $M^{A,B}=U(\slt\T \C[t]/t^n)\cdot (v_A\T v_B)\subset M^A\T M^B.$
Then $M^{A,B}\simeq M^C$ as $\slt\T \C[t]/t^n$ modules.
\end{prop}
\begin{proof}
Recall that in \cite{mi2} we constructed a fermionic realization of $M^A$
in the space $F^{\T (a_n-1)}$ -- the tensor power of the space of the
semi-infinite forms. Let us briefly describe this construction.

Consider the set of variables $\psi_n, \phi_m, n,m\in\Z$
with a relations:
$$[\psi_n,\psi_m]_+=[\phi_n,\phi_m]_+=[\psi_n,\phi_m]_+=0$$
(here $[a,b]_+=ab+ba$). Then the following vectors form a base of $F$:
$$\ld \psi_{N+1} \phi_{N+1} \psi_N \phi_N \psi_{i_1}\ld \psi_{i_k}
\phi_{j_1}\ld\phi_{j_l},\ N>i_1>\ld >i_k, N>j_1>\ld >j_l. $$
Note that operators $\psi_n, \phi_m$ act on $F$ by multiplication.
There are also a collection of operators $\psi^*_n, \phi^*_m$, acting on $F$
with a commutation relations
$$[\psi^*_n,\psi_m]_+=\delta_{n,-m},\
  [\phi^*_n,\phi_m]_+=\delta_{n,-m},\  [\phi^*_n,\phi^*_m]_+=0,
  [\psi^*_n,\psi^*_m]_+=0.$$
The Lie algebra $\slth$ acts on $F$ and
the action of $e_i, f_j\in\slth$ is given in terms of generating functions
\begin{equation*}
\psi(z)=\sum_{n\in\Z} \psi_n z^n,\ \phi(z)=\sum_{n\in\Z} \phi_n z^n,
\psi^*(z)=\sum_{n\in\Z} \psi^*_n z^n,\
\phi^*(z)=\sum_{n\in\Z} \phi_n z^n
\end{equation*}
in the following way:
$$e(z)=\sum_{n\in\Z} e_nz^n=\psi(z)\phi(z),\
f(z)=\sum_{n\in\Z} e_nz^n=\psi^*(z)\phi^*(z).$$
Note that $F$ is level $1$ module.

Introduce a notation for the extremal vectors:
$$v_{2N}=\ld \psi_{N+2}\phi_{N+2}\psi_{N+1}\phi_{N+1}\psi_N,\
  v_{2N+1}=\ld \psi_{N+2}\phi_{N+2}\psi_{N+1}\phi_{N+1}.$$
For $A\in\N^n$ denote $d_i=\#\{j:\ a_j=i\}$. In \cite{mi2} it was
proved that
$$M^A\simeq U(\slt\T\C[t])\cdot (v_{n-d_1}\T v_{n-d_1-d_2}\T\ld\T
v_{n-d_1-\ld-d_{a_n-1}})\subset F^{\T (a_n-1)}.$$
So we have an embedding $\imath: M^A\to F^{\T (a_n-1)}$. To
prove the proposition it is enough to mention that
$\imath (v_A)\T \imath (v_B) =\imath (v_C)$.
\end{proof}
\begin{cor}
\label{inbet}
For any $n>k\ge 1$ there exists a bundle $\pi_{n,k}:\shn\to \sh^{(k)}$ with
a fiber $\sh^{(n-k)}$.
\end{cor}
\begin{proof}
We use the induction on $n$.
Let $\shn$  be realized as $\sh_{(2,3,\ld,n+1)}$. From the previous lemma we
obtain an embedding for any $k<n$:
$$\shn\hk\sh^{(k)}\times
\sh_{(2,3,\ld,n+1-k,\lefteqn{\underbrace{\phantom{n+1-k,\ld n }}_k}
n+1-k,\ld,n+1-k)}.$$
Thus, we obtain a map $\pi_{n,k}: \shn\to\sh^{(k)}$ as a composition of the
embedding and projection on the first factor.
Note that the following diagram is commutative:
$$
\begin{CD}
\shn @>\pi_{n,k}>> \sh^{(k)}\\
@V\pi_{n,1}VV @VV\pi_{k,1}V\\
\sh^{(1)} @= \sh^{(1)}
\end{CD}
$$
Fix an arbitrary point $x\in\sh^{(k)}$. From the commutative diagram one can
see that $\pi_{n,k}^{-1}(x)$ is isomorphic to the preimage of the point
via the map of the fibers of $\pi_{n,1}$ and $\pi_{k,1}$, i.e.
the map between $\sh^{(n-1)}$ and $\sh^{(k-1)}$. But this map is
$\pi_{n-1,k-1}$. By the induction, we know that the fiber of this bundle
is $\sh^{(n-k)}$. Corollary is proved.
\end{proof}

\subsection{Identification of $\shn$ with the generalized flag manifold.\\}
Let $V_0=\C^2\T (\C[t]/t^n)$ with a natural action of
the Lie algebra $\slt\T(\C[t]/t^n)$, the Lie group
$\G=\Slt(\C[t]/t^n)$ and the operator $t$, acting by the multiplication.
Fix a basis $v_+, v_-$ of $\C^2$
with $h_0v_+=-v_+$, $h_0v_-=v_-$.
Define ${\bf F}_n$ as a variety of the following
sequences of the subspaces of $V_0$:
$$V_1\hookleftarrow V_2\hookleftarrow\ld
\hookleftarrow V_n$$
with a properties
$$tV_i\hk V_{i+1}\quad \text{ and }\quad \dim (V_i/V_{i+1})=1,\qquad
i=0,\ld,n-1.$$ Note that the group
$\Slt(\C[t]/t^n)$ acts on ${\bf F}_n$. Define the points $[P]$ and $[Q]$
in
${\bf F}_n$ as follows:
\begin{gather*}
[P]=V_1\hki \ld \hki V_n,\quad V_i=\bra v_+\T t^k, v_-\T t^l,\
0\le k\le n-1,\ i\le l \le n-1\ket;\\
[Q]=U_1\hki \ld \hki U_n,\quad U_i=\bra v_+\T t^k, v_-\T t^l,\
i\le k\le n-1,\ 0\le l\le n-1\ket.
\end{gather*}
Note that ${\bf F}_1\simeq \Pro^1.$
The following lemmas contain the list of the statements, which can be
checked directly.
\begin{lem}
\label{locrec1}
\begin{enumerate}
\item ${\bf F}_n$ is fibered over ${\bf F}_1\simeq \Pro^1$ with a fiber
${\bf F}_{n-1}$.
\item For any $n>k$ there exists a bundle ${\bf F}_n\to {\bf F}_k$ with a
fiber ${\bf F}_{n-k}$.
\end{enumerate}
\end{lem}

In the following lemma $w$ is a Weyl element from
$\Slt\hk \Slt(\C[t]/t^n)=\G$ and $\U_+$ and $\U_-$ are the subgroups of $\G$,
defined in $(\ref{unip})$.
\begin{lem}
\label{locrec2}
\begin{enumerate}
\item $w [P]=[Q]$.
\item $\U_+\cdot [P]\simeq \C^n\simeq \U_-\cdot [Q].$
\item $\G\cdot [P]=\U_+\cdot [P]\cup \U_-\cdot [Q].$
\item The orbit $\G\cdot [P]$ is fibered over $\Pro^1$ with a
fiber $\C^{n-1}$.
\item $\ov{\U_+\cdot [P]}={\bf F}_n=\ov{\U_-\cdot [Q]}$.
\end{enumerate}
\end{lem}

\begin{theorem}
There is a $\G$-isomorphism ${\bf F}_n\to \shn$, sending $[P]$ to $[v_A]$
(we fix some realization of $\shn$ as $\sh_A$).
\end{theorem}
\begin{proof}
Because of the lemmas $(\ref{locrec1}), (\ref{locrec2})$ ${\bf F}_n$
is determined by the
transition functions of the bundle ${\bf F}_n\to\Pro^1$.
The latter are determined by the transition functions of
$\G\cdot [P]\to\Pro^1$.
By the same
arguments as in the lemma $(\ref{bunlem})$ we obtain that these
functions coincide with the ones for $\shn$.
The theorem is proved.
\end{proof}

\subsection{Vector fields on $\shn$.\\ }
In this subsection we calculate the Lie algebra $\Vect (\shn)$ of the vector
fields on $\shn$.

It is clear that there is an embedding
$\slt\T(\C[t]/t^n)\hk \Vect (\shn).$ Note also that there are operators
$L_i, i=0,1\ld,n-2$, acting on $M^A$. Recall that in \cite{mi2} for any
$A$ we have constructed an integrable $\slth$ module with $M^A$ as
an $\sltc$ submodule. Moreover, for the operators $L_i$ from the Virasoro
algebra, acting on any integrable affine algebra module, we have
$$L_iv_A=0,\ i>0;\ \ [L_i,e_j]=-je_{i+j}.$$
Thus, the subalgebra of the Virasoro algebra, spanned by $L_i, i\ge 0$
is
acting on $M^A$ (for the convenience, we put $L_0v_A=0$; note also that
$L_{>n-2}$ acts on $M^A$ by $0$). In the
following lemma we prove that the exponents of the operators
$L_i, 0\le i\le n-2$
define a vector fields on $\shn$ and write an explicit formulas for the
corresponding vector fields on the big cell.
\begin{rem}
Note that using the identification of $\shn$ with the
generalized flag manifold one obtains a natural action of the vector
fields $L_i$ on $\shn$.
\end{rem}

\begin{lem}
\label{vf}
$1).$ Let $[x]\in\shn$. Then the vector
$$\frac{\partial \exp(L_i\veps)[x]}{\partial\veps}|_{\veps=0}\in
T_{[x]}\Pro (M^A)$$ is an element of $T_{[x]}{\shn}$. The corresponding
vector fields on a big cell
$$\left\{\exp \left(\sum_{i=0}^{n-1} e_ix_i\right),\ x_i\in\C\right\}$$
are given by the formula
\begin{equation}
\label{vect}
L_i=\sum_{j=1}^{n-i-1} jx_j\partial_{x_{i+j}},\qquad 0\le i\le n-2.
\end{equation}
$2).$ The restriction of the vector fields from $\slt\T(\C[t]/t^n)$ on a big
cell is
given by the formulas ($i=0,\ld,n-1$):
\begin{gather}
\label{e}
e_i=\partial_{x_i};\\
\label{h}
h_i=-2\sum_{j=0}^{n-i-1} x_j\partial_{x_{i+j}};\\
\label{f}
f_i=-\sum_{j=0}^{n-1-i} \left( \sum_{a+b=j} x_ax_b \right)
\partial_{x_{i+j}}.
\end{gather}
\end{lem}
\begin{proof}
In order to prove that the operators $L_i$ define a vector fields on $\shn$
it is enough to show that the latter is true on a big cell.
Consider an operator $\bar L_i$ acting on the space with a basis $e_j$ by the
formula
$\bar L_i (e_j)=-j e_{i+j}$. Then we obtain
\begin{equation}
\label{exp1}
\exp(L_i\veps)\left(\sum_{j=0}^{n-1} e_jx_j\right)\exp(-L_i\veps)=
\exp(\bar L_i\veps)\left(\sum_{j=0}^{n-1} e_jx_j\right)
\end{equation}
and thus
\begin{equation}
\label{exp2}
\exp(L_i\veps)\exp\left(\sum_{j=0}^{n-1} e_jx_j\right)\exp(-L_i\veps)=
\exp\left(\exp(\bar L_i\veps)\left(\sum_{j=0}^{n-1} e_jx_j\right)\right).
\end{equation}
The right hand side of $(\ref{exp1})$ is a series in $e_j$. Its exponent,
which is the right hand side of $(\ref{exp2})$, is a series in $e_j$ too.
Now apply both right and left hand sides of $(\ref{exp2})$ to the point
$[v_A]$. This gives us that the operators $L_i$ really define a vector fields
on the  big cell (the only thing to use is $L_i v_A=0$).
Note also that $(\ref{vect})$ is a trivial consequence of $(\ref{exp2})$.
Thus, the first part of our lemma is proved.

The only thing to prove in the second part, is the representation of
$f_i$ as a vector field.
Let $x(z)=x_0+zx_1+\ld +z^{n-1}x_{n-1}, y(z)=y_0+zy_1+\ld +z^{n-1} y_{n-1}$
be such elements of $\C[z]/z^n$ that $x(z)y(z)=1$.
Recall, (see lemma $(\ref{bunlem})$) that we have
the following equality:
$$\exp\left(\sum_{i=0}^{n-1} e_i x_i\right)[v_A]=
\exp\left(\sum_{i=0}^{n-1} f_i y_i\right)[u_A].
$$
Hence we obtain
\begin{multline}
\exp(f_i\veps)\exp\left(\sum_{k=0}^{n-1} e_kx_k\right)[v_A]=
\exp(f_i\veps)\exp\left(\sum_{k=0}^{n-1} f_ky_k\right)[u_A]=\\
=\exp\left(\sum_{k\ne i} f_ky_k +f_i(y_i+\veps)\right)[u_A]=
\exp\left(\sum_{k=0}^{n-1} e_k x'_k\right)[v_A].
\end{multline}
Here for $x'(z)=\sum_{k=0}^{n-1} x'_kz^k$ we have
$$x'(z)=\frac{1}{\veps z^i+x(z)^{-1}}.$$
This gives $(\ref{f})$. To obtain $(\ref{h})$ one must commute the
expressions for $e_i$ and $f_j$.
Lemma is proved.
\end{proof}

Now our goal is to prove that  $\bra e_i, h_i, L_i, f_i\ket =\Vect(\shn)$.
To do that, it is enough to show that $\dim \Vect(\shn)=4n-1$
(the linear independence is a trivial consequence from
$(\ref{vect})$, $(\ref{e})$, $(\ref{h})$, and $(\ref{f})$).

Let $T_n$ be a tangent bundle on $\shn$. Recall that there is a bundle
$\pi_{n,1}: \shn\to \Pro^1$ with a fiber $\sh^{(n-1)}$. Thus we  obtain a
surjection
$\Vect(\shn)\sur \Vect(\Pro^1)\simeq \slt$. Its kernel consists of the vector
fields on $\shn$ which are tangent to the fibers of $\pi_{n,1}$. Denote
by $T'_n$ the bundle on $\shn$, whose fiber at the point $x$ is a tangent
space to the fiber $\pi_{n,1}^{-1}(\pi_{n,1}x)$. We need to prove that
$\dim H^0 (T'_n)=4n-4.$ In order to do that, consider a bundle $\E_n$ on
$\Pro^1$, whose fiber at the point $z$ is the space of vector fields on
$\pi_{n,1}^{-1}(z)$ (thus, $\E_n=(\pi_{n,1})_* T_n'$). It is obvious that
$\dim H^0 (T'_n)=\dim H^0 (\E_n)$. To prove that the latter equals to
$4n-4$ we will compute the transition functions of $\E_n$. The idea is
to use the induction on $n$, i.e. the knowledge of the Lie algebra of
the vector fields on a fiber of $\pi_{n,1}$.  As a result we will prove
that:
\begin{gather*}
\E_n=\Ob(2)\oplus \Ob(1)^{\oplus (n-1)}\oplus \Ob(0)^{\oplus (2n-5)}\oplus
\Ob(-1)^{\oplus (n-1)}\oplus \Ob(-2),\qquad n>2,\\
\E_2=\Ob(2)\oplus\Ob(0)\oplus \Ob(-2).
\end{gather*}

Recall that we have two big cells in $\shn$:
$$U_x=\exp\left(\sum_{i=0}^{n-1} x_ie_i\right)[v_A]\simeq \C^n\simeq
\exp\left(\sum_{i=0}^{n-1} y_if_i\right)[u_A]=U_y.$$
Let $e_{x,i}, f_{x,i}, h_{x,i}, L_{x,i}$ be the vector fields on $U_x$, which
are the restrictions of the vector fields $e_i, f_i, h_i, L_i$
correspondingly. Denote also by $e_{y,i}, f_{y,i}, h_{y,i}, L_{y,i}$ the
vector fields on $U_y$, which are the restrictions of the vector fields
$f_i, e_i, -h_i, L_i$ correspondingly (this notations are convenient for us,
because of the simplification of some formulas).
Surely, on $U_x\cap U_y$ we have
\begin{equation}
\label{coor}
e_{x,i}=f_{y,i},\quad f_{x,i}=e_{y,i},\quad h_{x,i}=-h_{y,i}, \quad
L_{x,i}=L_{y,i}.
\end{equation}
We can write an explicit formulas for the above vector fields on $U_y$:
\begin{gather}
e_{y,i}=\partial_{y_i},\qquad
h_{y,i}=-2\sum_{j=0}^{n-i-1} y_j\partial_{y_{i+j}}, \\
f_{y,i}=-\sum_{j=0}^{n-i-1} \left(\sum_{a+b=j} y_ay_b\right)\partial_{y_{i+j}},
\qquad L_{y,i}=\sum_{j=1}^{n-1-i} jy_j \partial_{y_{i+j}}.
\end{gather}

Suppose that we have already proved that $e_i, h_i, f_i (0\le i\le n-1)$ and
$L_i (0\le i\le n-2)$ form a base of $\Vect(\sh^{(n-1)})$.
Then we can trivialize the bundle $\E_n$ on $\pi_{n,1}U_x$, choosing the
following base of the vector fields, tangent to the fibers of $\pi_{n,1}$:
\begin{gather}
\label{vecx}
{e'}^x_i=\pa_{x_i},;\qquad
{h'}^x_i=-2\sum_{j=1}^{n-i} x_j\pa_{x_{i+j-1}};\\
{f'}^x_i=
-\sum_{j=1}^{n-i}\left(\sum_{\al+\be=j+1;\ \al,\be\ge 1}
x_\al x_\be\right) \pa_{x_{i+j-1}}, \qquad  i=1,\ld,n-1;\notag\\
{L'}^x_i=\sum_{j=1}^{n-i-1} jx_{j+1}\pa_{x_{i+j}},\ i=1,\ld,n-2.\notag
\end{gather}
We also trivialize $\E_n$ on $\pi_{n,1} U_y$, fixing the analogous basis:
\begin{gather}
\label{vecy}
{e'}^y_i=\pa_{y_i},;\qquad
{h'}^y_i=-2\sum_{j=1}^{n-i} y_j\pa_{y_{i+j-1}};\\
{f'}^y_i=
-\sum_{j=1}^{n-i}\left(\sum_{\al+\be=j+1;\ \al,\be\ge 1}
y_\al y_\be\right) \pa_{y_{i+j-1}}, \qquad  i=1,\ld,n-1;\notag\\
{L'}^y_i=\sum_{j=1}^{n-i-1} jy_{j+1}\pa_{y_{i+j}},\ i=2,\ld,n-1.\notag
\end{gather}
To obtain the transition functions of $\E_n$ one must rewrite the
vector fields $(\ref{vecx})$ on
$U_x\cap U_y$ via the vector fields from
$(\ref{vecy})$.

\begin{lem}
\label{n'}
We have the following equalities of the vector fields on $U_x$:
\begin{gather*}
e_{x,i}={e'}^x_i,\qquad h_{x,i}={h'}^x_{i+1}-2x_0{e'}^x_i,\\
L_{x,i}={L'}^x_{i+1}-{h'}^x_i/2,\qquad
f_{x,i}={f'}^x_{i+2}+x_0{h'}^x_{i+1}-x_0^2 {e'}^x_i.
\end{gather*}
The analogous formulas are true on $U_y$:
\begin{gather*}
e_{y,i}={e'}^y_i,\qquad h_{y,i}={h'}^y_{i+1}-2y_0{e'}^y_i,\\
L_{y,i}={L'}^y_{i+1}-{h'}^y_i/2,\qquad
f_{y,i}={f'}^y_{i+2}+y_0{h'}^y_{i+1}-y_0^2 {e'}^y_i.
\end{gather*}
\end{lem}
\begin{proof}
It is a consequence from the lemma $(\ref{vf})$ and the definitions
$(\ref{vecx})$, $(\ref{vecy})$.
\end{proof}

\begin{lem}
We have the following equalities (note that $y_0^{-1}=x_0$):
\begin{gather}
\label{xy}
{e'}^x_i=y_0^2{e'}^y_i+y_0{h'}^y_{i+1}+{f'}^y_{i+2},\qquad
{h'}^x_i={h'}^y_i+2y_0^{-1}{f'}^y_{i+1},\\
{L'}^x_i={L'}^y_i+y_0^{-1}{f'}^y_{i+1},\qquad {f'}^x_i=-y_0^{-2}{f'}^y_i.
\notag
\end{gather}
\end{lem}
\begin{proof}
Using the formulas from the lemma $(\ref{n'})$ we can rewrite our
fields, written in the $x_i$-coordinates via the $y_i$-coordinates:
\begin{multline}
{e'}^x_i=e_{x,i}=f_{y,i}={f'}^y_{i+2}+y_0{h'}^y_{i+1}-y_0^2{e'}^y_i;\\
\shoveleft {{h'}^x_i=h_{x,i-1}+2x_0e_{x,i-1}=-h_{y,i-1}+2y_0^{-1}f_{y,i-1}=}\\
=-{h'}^y_i+2y_0{e'}^y_{i-1}+2y_0^{-1}({f'}^y_{i+1}+y_0{h'}^y_i-
y_0^2{e'}^y_{i-1})={h'}^y_i+2y_0^{-1}{f'}^y_{i+1};\\
\shoveleft
{{L'}^x_i=L_{x,i-1}+{h'}^x_{i-1}/2=L_{y,i-1}+{h'}^y_{i-1}/2+
y_0^{-1}{f'}^y_i=}\\
={L'}^y_i-{h'}^y_{i-1}/2+{h'}^y_{i-1}/2+y_0^{-1}{f'}^y_i=
{L'}^y_i+y_0^{-1}{f'}^y_i;\\
\shoveleft {{f'}^x_i=f_{x,i-2}-x_0h_{x,i-2}-x_0^2e_{x,i-2}=
e_{y,i-2}+y_0^{-1}h_{y,i-2}-y_0^{-2}f_{y,i-2}=}\\
={e'}^y_{i-2}+y_0^{-1}({h'}^y_{i-1}-2y_0{e'}^y_{i-2})-
y_0^{-2}({f'}^y_i+y_0{h'}^y_{i-1}-y_0^2 {e'}^y_{i-2})=-y_0^{-2}{f'}^y_i.
\end{multline}
This gives us the formulas $(\ref{xy})$.
\end{proof}

Now we are ready to prove the main theorem.
\begin{theorem}
\begin{gather*}
\E_n=\Ob(2)\oplus \Ob(1)^{\oplus (n-1)}\oplus \Ob(0)^{\oplus (2n-5)}\oplus
\Ob(-1)^{\oplus (n-1)}\oplus \Ob(-2),\qquad n>2,\\
\E_2=\Ob(2)\oplus\Ob(0)\oplus \Ob(-2).
\end{gather*}
\end{theorem}
\begin{proof}
Let $n>2$. The transition matrix of $\E_n$
has the following form:
{\footnotesize{
$$\begin{array}{c|cccccccccccccccccc}
 &\!\!\! {e'}_1^x\!\!\! &\!\!\! \ldots \!\!\! &\!\!\! {e'}_{n-3}^x \!\!\!&
 \!\!\!
 {e'}_{n-2}^x\!\!\!
 &\!\!\! {e'}_{n-1}^x\!\! & \!\! {h'}_1^x \!\!\! & \!\!\! {h'}_2^x \!\!\! &
 \! \!\! \ldots\!\! \! &\!\! \!
 {h'}_{n-2}^x\!\!  \!
 &\!\! \! {h'}_{n-1}^x\!\!  & \!\! {L'}_1^x\!\! \! & \!\!\! \ldots
 \!\! \! &\!\! \!
 {L'}_{n-2}^x\!\!
 &\!\!  {f'}_1^x\!\! \!  &\!\! \! {f'}_2^x\!\! \! & \!\! \! {f'}_3^x \!\!\!
 &
 \!\!\!  \ldots\!\! \! &\!\! \! {f'}_{n-1}^x\\
\hline

\!\!\! {e'}_1^y\!\!\!
&\!\!\! -y_0^2\!\!\! &\!\!\! \ldots \!\!\! &\!\!\! 0 \!\!\!&
 \!\!\!
 0\!\!\!
 &\!\!\! 0\!\! & \!\! 0 \!\!\! & \!\!\! 0 \!\!\! &
 \! \!\! \ldots\!\! \! &\!\! \!
 0\!\!  \!
 &\!\! \! 0\!\!  & \!\! 0\!\! \! & \!\!\! \ldots
 \!\! \! &\!\! \!
 0\!\!
 &\!\!  0\!\! \!  &\!\! \! 0\!\! \! & \!\! \! 0 \!\!\!
 &
 \!\!\!  \ldots\!\! \! &\!\! \! 0\\

\!\!\! \vdots \!\!\!& \vdots &\!\! \ddots \!\!\!&\!\! \vdots\!\!\! &\!\!\! \vdots\!\!
&\!\! \vdots
&\!\! \vdots \!\!&
\!\!\vdots\! & \!\! \ddots\!\!
&
 \vdots\!\!\! &\!\!\! \vdots \!\!\!&\!\!\! \vdots\! &\!\!\! \ddots\!\!\! &
 \!\vdots\!
 &\! \vdots\! &
\!\vdots \!&\! \vdots & \!\! \ddots\!\! &
 \vdots\\

\!\!\! {e'}_{n-3}^y \!\!\!
&\!\!\! 0\!\!\! &\!\!\! \ldots \!\!\! &\!\!\! -\!y_0^2 \!\!\!&
 \!\!\!
 0\!\!\!
 &\!\!\! 0\!\! & \!\! 0 \!\!\! & \!\!\! 0 \!\!\! &
 \! \!\! \ldots\!\! \! &\!\! \!
 0\!\!  \!
 &\!\! \! 0\!\!  & \!\! 0\!\! \! & \!\!\! \ldots
 \!\! \! &\!\! \!
 0\!\!
 &\!\!  0\!\! \!  &\!\! \! 0\!\! \! & \!\! \! 0 \!\!\!
 &
 \!\!\!  \ldots\!\! \! &\!\! \! 0\\

\!\!\! {e'}_{n-2}^y \!\!\!
&\!\!\! 0\!\!\! &\!\!\! \ldots \!\!\! &\!\!\! 0 \!\!\!&
 \!\!\!
 -y_0^2\!\!\!
 &\!\!\! 0\!\! & \!\! 0 \!\!\! & \!\!\! 0 \!\!\! &
 \! \!\! \ldots\!\! \! &\!\! \!
 0\!\!  \!
 &\!\! \! 0\!\!  & \!\! 0\!\! \! & \!\!\! \ldots
 \!\! \! &\!\! \!
 0\!\!
 &\!\!  0\!\! \!  &\!\! \! 0\!\! \! & \!\! \! 0 \!\!\!
 &
 \!\!\!  \ldots\!\! \! &\!\! \! 0\\

\!\!\! {e'}_{n-1}^y \!\!\!
&\!\!\! 0\!\!\! &\!\!\! \ldots \!\!\! &\!\!\! 0 \!\!\!&
 \!\!\!
 0\!\!\!
 &\!\!\! -y_0^2\!\! & \!\! 0 \!\!\! & \!\!\! 0 \!\!\! &
 \! \!\! \ldots\!\! \! &\!\! \!
 0\!\!  \!
 &\!\! \! 0\!\!  & \!\! 0\!\! \! & \!\!\! \ldots
 \!\! \! &\!\! \!
 0\!\!
 &\!\!  0\!\! \!  &\!\! \! 0\!\! \! & \!\! \! 0 \!\!\!
 &
 \!\!\!  \ldots\!\! \! &\!\! \! 0\\

\!\!\! {h'}_1^y  \!\!\!
&\!\!\! 0\!\!\! &\!\!\! \ldots \!\!\! &\!\!\! 0 \!\!\!&
 \!\!\!
 0\!\!\!
 &\!\!\! 0\!\! & \!\! 1 \!\!\! & \!\!\! 0 \!\!\! &
 \! \!\! \ldots\!\! \! &\!\! \!
 0\!\!  \!
 &\!\! \! 0\!\!  & \!\! 0\!\! \! & \!\!\! \ldots
 \!\! \! &\!\! \!
 0\!\!
 &\!\!  0\!\! \!  &\!\! \! 0\!\! \! & \!\! \! 0 \!\!\!
 &
 \!\!\!  \ldots\!\! \! &\!\! \! 0\\

\!\!\! {h'}_2^y  \!\!\!
&\!\!\! y_0\!\!\! &\!\!\! \ldots \!\!\! &\!\!\! 0 \!\!\!&
 \!\!\!
 0\!\!\!
 &\!\!\! 0\!\! & \!\! 0 \!\!\! & \!\!\! 1 \!\!\! &
 \! \!\! \ldots\!\! \! &\!\! \!
 0\!\!  \!
 &\!\! \! 0\!\!  & \!\! 0\!\! \! & \!\!\! \ldots
 \!\! \! &\!\! \!
 0\!\!
 &\!\!  0\!\! \!  &\!\! \! 0\!\! \! & \!\! \! 0 \!\!\!
 &
 \!\!\!  \ldots\!\! \! &\!\! \! 0\\

\!\!\! \vdots \!\!\!& \vdots &\!\! \ddots \!\!\!&\!\! \vdots\!\!\! &\!\!\!
\vdots\!\!
&\!\! \vdots
&\!\! \vdots \!\!&
\!\!\vdots\! & \!\! \ddots\!\!
&
 \vdots\!\!\! &\!\!\! \vdots \!\!\!&\!\!\! \vdots\! &\!\!\! \ddots\!\!\! &
 \!\vdots\!
 &\! \vdots\! &
\!\vdots \!&\! \vdots & \!\! \ddots\!\! &
 \vdots\\

\!\!\! {h'}_{n-2}^y \!\!\!
&\!\!\! 0\!\!\! &\!\!\! \ldots \!\!\! &\!\!\! y_0 \!\!\!&
 \!\!\!
 0\!\!\!
 &\!\!\! 0\!\! & \!\! 0 \!\!\! & \!\!\! 0 \!\!\! &
 \! \!\! \ldots\!\! \! &\!\! \!
 1\!\!  \!
 &\!\! \! 0\!\!  & \!\! 0\!\! \! & \!\!\! \ldots
 \!\! \! &\!\! \!
 0\!\!
 &\!\!  0\!\! \!  &\!\! \! 0\!\! \! & \!\! \! 0 \!\!\!
 &
 \!\!\!  \ldots\!\! \! &\!\! \! 0\\

\!\!\! {h'}_{n-1}^y \!\!\!
&\!\!\! 0\!\!\! &\!\!\! \ldots \!\!\! &\!\!\! 0 \!\!\!&
 \!\!\!
 y_0\!\!\!
 &\!\!\! 0\!\! & \!\! 0 \!\!\! & \!\!\! 0 \!\!\! &
 \! \!\! \ldots\!\! \! &\!\! \!
 0\!\!  \!
 &\!\! \! 1\!\!  & \!\! 0\!\! \! & \!\!\! \ldots
 \!\! \! &\!\! \!
 0\!\!
 &\!\!  0\!\! \!  &\!\! \! 0\!\! \! & \!\! \! 0 \!\!\!
 &
 \!\!\!  \ldots\!\! \! &\!\! \! 0\\

\!\!\! {L'}_1^y  \!\!\!
&\!\!\! 0\!\!\! &\!\!\! \ldots \!\!\! &\!\!\! 0 \!\!\!&
 \!\!\!
 0\!\!\!
 &\!\!\! 0\!\! & \!\! 0 \!\!\! & \!\!\! 0 \!\!\! &
 \! \!\! \ldots\!\! \! &\!\! \!
 0\!\!  \!
 &\!\! \! 0\!\!  & \!\! 1\!\! \! & \!\!\! \ldots
 \!\! \! &\!\! \!
 0\!\!
 &\!\!  0\!\! \!  &\!\! \! 0\!\! \! & \!\! \! 0 \!\!\!
 &
 \!\!\!  \ldots\!\! \! &\!\! \! 0\\

\!\!\! \vdots \!\!\!& \vdots &\!\! \ddots \!\!\!&\!\! \vdots\!\!\! &\!\!\! \vdots\!\!
&\!\! \vdots
&\!\! \vdots \!\!&
\!\!\vdots\! & \!\! \ddots\!\!
&
 \vdots\!\!\! &\!\!\! \vdots \!\!\!&\!\!\! \vdots\! &\!\!\! \ddots\!\!\! &
 \!\vdots\!
 &\! \vdots\! &
\!\vdots \!&\! \vdots & \!\! \ddots\!\! &
 \vdots\\

\!\!\! {L'}_{n-2}^y \!\!\!
&\!\!\! 0\!\!\! &\!\!\! \ldots \!\!\! &\!\!\! 0 \!\!\!&
 \!\!\!
 0\!\!\!
 &\!\!\! 0\!\! & \!\! 0 \!\!\! & \!\!\! 0 \!\!\! &
 \! \!\! \ldots\!\! \! &\!\! \!
 0\!\!  \!
 &\!\! \! 0\!\!  & \!\! 0\!\! \! & \!\!\! \ldots
 \!\! \! &\!\! \!
 1\!\!
 &\!\!  0\!\! \!  &\!\! \! 0\!\! \! & \!\! \! 0 \!\!\!
 &
 \!\!\!  \ldots\!\! \! &\!\! \! 0\\

\!\!\! {f'}_1^y  \!\!\!
&\!\!\! 0\!\!\! &\!\!\! \ldots \!\!\! &\!\!\! 0 \!\!\!&
 \!\!\!
 0\!\!\!
 &\!\!\! 0\!\! & \!\! 0 \!\!\! & \!\!\! 0 \!\!\! &
 \! \!\! \ldots\!\! \! &\!\! \!
 0\!\!  \!
 &\!\! \! 0\!\!  & \!\! 0\!\! \! & \!\!\! \ldots
 \!\! \! &\!\! \!
 0\!\!
 &\!\!  \frac{-1}{y_0^2}\!\! \!  &\!\! \! 0\!\! \! & \!\! \! 0 \!\!\!
 &
 \!\!\!  \ldots\!\! \! &\!\! \! 0\\

\!\!\! {f'}_2^y  \!\!\!
&\!\!\! 0\!\!\! &\!\!\! \ldots \!\!\! &\!\!\! 0 \!\!\!&
 \!\!\!
 0\!\!\!
 &\!\!\! 0\!\! & \!\! \frac{-2}{y_0} \!\!\! & \!\!\! 0 \!\!\! &
 \! \!\! \ldots\!\! \! &\!\! \!
 0\!\!  \!
 &\!\! \! 0\!\!  & \!\! \frac{1}{y_0}\!\! \! & \!\!\! \ldots
 \!\! \! &\!\! \!
 0\!\!
 &\!\!  0\!\! \!  &\!\! \! \frac{-1}{y_0^2}\!\! \! & \!\! \! 0 \!\!\!
 &
 \!\!\!  \ldots\!\! \! &\!\! \! 0\\

\!\!\! {f'}_3^y  \!\!\!
&\!\!\! 1\!\!\! &\!\!\! \ldots \!\!\! &\!\!\! 0 \!\!\!&
 \!\!\!
 0\!\!\!
 &\!\!\! 0\!\! & \!\! 0 \!\!\! & \!\!\! \frac{-2}{y_0} \!\!\! &
 \! \!\! \ldots\!\! \! &\!\! \!
 0\!\!  \!
 &\!\! \! 0\!\!  & \!\! \frac{1}{y_0}\!\! \! & \!\!\! \ldots
 \!\! \! &\!\! \!
 0\!\!
 &\!\!  0\!\! \!  &\!\! \! 0\!\! \! & \!\! \! \frac{-1}{y_0^2} \!\!\!
 &
 \!\!\!  \ldots\!\! \! &\!\! \! 0\\

\!\!\! \vdots \!\!\!& \vdots &\!\! \ddots \!\!\!&\!\! \vdots\!\!\! &\!\!\! \vdots\!\!
&\!\! \vdots
&\!\! \vdots \!\!&
\!\!\vdots\! & \!\! \ddots\!\!
&
 \vdots\!\!\! &\!\!\! \vdots \!\!\!&\!\!\! \vdots\! &\!\!\! \ddots\!\!\! &
 \!\vdots\!
 &\! \vdots\! &
\!\vdots \!&\! \vdots & \!\! \ddots\!\! &
 \vdots\\

\!\!\! {f'}_{n-1}^y \!\!\!
&\!\!\! 0\!\!\! &\!\!\! \ldots \!\!\! &\!\!\! 1 \!\!\!&
 \!\!\!
 0\!\!\!
 &\!\!\! 0\!\! & \!\! 0 \!\!\! & \!\!\! 0 \!\!\! &
 \! \!\! \ldots\!\! \! &\!\! \!
 \frac{-2}{y_0}\!\!  \!
 &\!\! \! 0\!\!  & \!\! 0\!\! \! & \!\!\! \ldots
 \!\! \! &\!\! \!
\frac{1}{y_0}\!\!
 &\!\!  0\!\! \!  &\!\! \! 0\!\! \! & \!\! \! 0 \!\!\!
 &
 \!\!\!  \ldots\!\! \! &\!\! \! \frac{-1}{y_0^2}\\
\end{array}$$}}
Now to identify the bundle $\E_n$ for $n>2$ one must diagonalize this matrix
using two following operations:
\begin{enumerate}
\item To add to some column another one, multiplied by the polynomial
in the variable $y_0^{-1}$.
\item To add to some row another one, multiplied by the polynomial in the
variable $y_0$.
\end{enumerate}
Recall that $y_0$ and $x_0=y_0^{-1}$ are the coordinates on the two charts of
$\Pro^1$: $\pi_{n,1}U_y$ and $\pi_{n,1}U_x$  respectively.
Thus the first operation is a change of the basis in $\Vect(U_x)$ and the
second one
is the change of the basis in $\Vect(U_y)$. It is easy to show that using
the above operations one can obtain the diagonal matrix:
$$\diag (y_0^2,\underbrace{y_0,\ld, y_0}_{n-1},\underbrace{1,\ld,1}_{2n-5},
\underbrace{y_0^{-1},\ld, y_0^{-1}}_{n-1}, y_0^{-2}).$$
This gives us our lemma in the case $n>2$. The case $n=2$ can be considered
in the same way.
\end{proof}
\begin{cor}
$\dim H^0(T'_n)=4n-4.$
\end{cor}
\begin{theorem}
\label{vectf}
The Lie algebra $\Vect(\shn)$ has a following basis:
$e_i, h_i, f_i$, $i=0,\ld,n-1$ and $L_i, i=0,\ld,n-2$.
\end{theorem}
\begin{proof}
As we have already mentioned, the only thing to prove is
$\dim \Vect(\shn)=4n-1$. But $\dim \Vect(\shn)=
\dim H^0 (T'_n)+\dim \Vect (\Pro^1)=4n-4+3=4n-1.$ Theorem is proved.
\end{proof}

\section{Back to the fusion.}
For the successive study of the Schubert varieties we need some additional
information
about the structure of the fusion products. Namely, we need to study the
special submodules, mentioned in \cite{mi2}.

\subsection {Preliminaries.}
Recall that $M^A$, $A=(a_1,\ld,a_n)\in \N^n$,
is $\sltc$ module and $\dim M^A=\prod_{i=1}^n a_i$.
Moreover, $M^A$ is cyclic $e\T\C[t]$ module
with a cyclic vector $v_A$ and defining relations
\begin{equation}
\label{power}
e_{(n)}(z)^i\div z^{\sum_{j=1}^n (i+1-a_j)_+},\ i=1,2,\ld
\end{equation}
($e_{(n)}(z)=\sum_{i=0}^{n-1} e_i z^{n-1-i}$, $a_+=\max(a,0)$)
(see \cite{mi1}).
The latter means that $M^A\simeq \Cn/I_A$ and the ideal $I_A$ is generated
by the elements
\begin{equation}
\sum_{i_1+\ld +i_k=s} e_{n-1-i_1}e_{n-1-i_2}\ld e_{n-1-i_k},\ k\ge 1,\
s<\sum_{j=1}^n (k+1-a_j)_+ .
\end{equation}

Let $a_1\le\ld\le a_n$.
In \cite{mi2} we constructed the following exact sequence
of $\slt\T \C[t]$ modules:
$$
0\to M^{(a_2-a_1+1,a_3,\ld,a_n)} \to M^A\to M^{(a_1-1,a_2+1,a_3,\ld,a_n)}
\to 0
$$
(note that $M^{(1,a_2,\ld,a_n)}\simeq M^{(a_2,\ld,a_n)}$).
One can generalize this construction in the following way.
For any $i,j$ with $1\le i<j\le n$ define
$$A_{i,j}=
(a_1,\ld,a_{i-1},a_i-1,a_{i+1},\ld,a_{j-1}, a_j+1,a_{j+1},\ld,a_n).$$
\begin{lem}
For any $1\le i<j\le n$ there exists an $\sltc$ module $S_{i,j}(A)\hk M^A$,
with a property
$M^A/S_{i,j}(A)\simeq M^{A_{i,j}}.$
\end{lem}
\begin{proof}
Note that the condition $(\ref{power})$ for the set
$A$ is weaker, than  for the set $A_{i,j}$.
Hence, there is a surjection
$\al_{i,j}(A): M^A\sur M^{A_{i,j}}$. We define $S_{i,j}(A)$ as its kernel.
Note that
$\al_{i,j}(A)$ is $\slt\T \C[t]$ homomorphism. In fact, it is obviously
$\Cn$ homomorphism. In addition,
$$f_i v_A=f_i v_{A_{i,j}}=0,\quad h_{>0} v_A= h_{>0} v_{A_{i,j}}=0$$
and the operator $h_0$ multiplies the vectors $v_A$ and $v_{A_{i,j}}$ by the
same constant. Thus, $S_{i,j}(A)$ is really the $\sltc$ module.
\end{proof}
\begin{rem}
\label{stop}
It was shown in \cite{mi2} that
\begin{enumerate}
\item $S_{1,2}(A)\simeq M^{(a_2-a_1+1,a_3,\ld,a_n)}$.
\item If $a_i=a_{i+1}$, then
$S_{i,i+1}(A)\simeq M^{(a_1,\ld,a_{i-1},a_{i+2},\ld,a_n)}.$
\end{enumerate}
\end{rem}
The rest of the section is devoted to the study of the modules $S_{i,i+1}(A)$.

\subsection{First description of $S_{i,i+1}(A)$.}
We construct a filtration on $S_{i,i+1}(A)$ with the quotients being the
fusion products.
Introduce a notation
$$A_i=(a_1,\ld, a_{i-2},a_{i-1}-a_i+a_{i+1}, a_{i+2},\ld,a_n), i>1.$$
\begin{lem}
We have an embedding $M^{A_i}\hk S_{i,i+1}(A)$.
\end{lem}
\begin{proof}
Denote $d_i=\#\{\al:\ a_\al=i\}, i=1,\ld,a_n.$  Let
$[e_{(n)}(z)^i]_j$ be a coefficient in $e_{(n)}(z)^i$ by the term $z^j$.
It was proved in \cite{mi2} that
the $\slt\T(\C[t]/t^{n-2})$ submodule in $M^A$ generated from the vector
\begin{equation}
\label{vec}
\left[e_{n}(z)^{a_i-1}\right]_{\sum_{j=1}^{a_i-1} (a_i-j)d_j} v_A
\end{equation}
is isomorphic to $M^{A_i}$.
One can show that the vector $(\ref{vec})$ is the element of the kernel of
surjection $M^A\sur M^{A_{i,i+1}}$.
Thus $M^{A_i}\hk S_{i,i+1}(A)$.
\end{proof}

Introduce a notation:
$N_A(k)=\sum_{j=1}^n (k+1-a_j)_+$.
In this notations the defining relations of $M^A$ read as
\begin{equation}
\label{N_A}
e_{(n)}(z)^k\div z^{N_A(k)},\quad k=1,2,\ld.
\end{equation}
We want to compare these relations in the case
of $M^A$ and $M^{A_{i,i+1}}$.
\begin{lem}
There are two cases:
\begin{enumerate}
\item  $a_i-1\le k\le a_{i+1}-1$. Then $N_A(k)=N_{A_{i,i+1}}(k)-1$.
\item  $k<a_i-1$ or $k>a_{i+1}-1$. Then $N_A(k)=N_{A_{i,i+1}}(k)$.
\end{enumerate}
\end{lem}
Recall that $S_{i,i+1}(A)$ is the kernel of the map $M^A\sur M^{A_{i,i+1}}$.
Hence, because of the formula $(\ref{N_A})$, we obtain that
$S_{i,i+1}(A)$ is generated by the
action of $\Cn$ from $a_{i+1}-a_i+1$ vectors, namely from
\begin{equation}
\label{genvec}
w_j=[e_{(n)}(z)^j]_{N_A(j)}v_A,\qquad j=a_i-1,\ld, a_{i+1}-1.
\end{equation}
For example, $\Cn\cdot w_{a_i-1}\simeq M^{A_i}\hk S_{i,i+1}(A)$.

\begin{lem}
\label{j-i}
Let $1\le i<j\le n$. Then $S_{i,j}(A)=S_{i,i+1}(A)+\ld+ S_{j-1,j}(A)$ (the
non-direct sum).
\end{lem}
\begin{proof}
It follows from the formula $(\ref{genvec})$.
\end{proof}

\begin{prop}
\label{filt}
Let $a_i\ne a_{i+1}, i\ne 1, a_{i-1}>1$. Then
\begin{multline}
\label{rec}
S_{i,i+1}(A)/M^{A_i}\simeq\\
\simeq S_{i,i+1}(a_1,\ld,a_{i-2},a_{i-1}-1,a_i+1,a_{i+1},\ld,a_n)=
S_{i,i+1}(A_{i-1,i}).
\end{multline}
\end{prop}
\begin{proof}
First, one can show that the dimensions of the both sides coincide. In fact,
by the definition
\begin{equation}
\label{first}
\dim S_{i,i+1}(A)=\dim M^A-\dim M^{A_{i,i+1}}=
\left(\prod_{j\ne i,i+1} a_j\right)\times (a_{i+1}-a_i+1),
\end{equation}
\begin{multline}
\label{sec}
\dim S_{i,i+1}(A_{i-1,i})=\dim M^{A_{i-1,i}}-
\dim M^{A_{i-1,i+1}}=\\
=\left(\prod_{j\ne i-1,i,i+1} a_j\right)\times (a_{i-1}-1)(a_{i+1}-a_i).
\end{multline}
It is easy to check that the difference between $(\ref{first})$ and
$(\ref{sec})$ equals to $\dim M^{A_i}$.

Consider the following mappings:
$$
M^{A_i} \stackrel{h}{\hk}
S_{i,i+1}(A)\stackrel{g}{\hk}
M^A\stackrel{\al_{i-1,i}(A)}{\longrightarrow} M^{A_{i-1,i}}
\stackrel{f}{\hookleftarrow} S_{i,i+1}(A_{i-1,i}).
$$
We will prove that
$$(\al_{i-1,i}(A)\circ g) (S_{i,i+1}(A))\subset f (S_{i,i+1}(A_{i-1,i}))$$
and $\ker (\al_{i-1,i}(A)\circ g)=h (M^{A_i})$. That will be enough for the
proof of the
proposition, because of the equality of the dimensions of the right
and left hand sides of $(\ref{rec})$.

We have the following commutative diagram:
$$
\begin{CD}
M^A @>\al_{i,i+1}(A)>> M^{A_{i,i+1}}\\
@V\al_{i-1,i}(A)VV @ VV\al_{i-1,i}(A_{i,i+1})V\\
M^{A_{i-1,i}} @>\al_{i,i+1}(A_{i-1,i})>> M^{A_{i-1,i+1}}
\end{CD}
$$
Note that $$S_{i,i+1}(A)=\ker \al_{i,i+1}(A)\quad \text{ and }\quad
S_{i,i+1}(A_{i-1,i})=\ker \al_{i,i+1} (A_{i-1,i}).$$
Hence,
$(\al_{i-1,i}(A)\circ g) (S_{i,i+1}(A))\subset f (S_{i,i+1}(A_{i-1,i})).$

Now, let us prove that $\ker \al_{i-1,i}(A)\circ g=h (M^{A_i})$.
Note that $\ker \al_{i-1,i}(A)=S_{i-1,i}(A)$.
Thus it is enough to show that
$$M^{A_i}=S_{i-1,i}(A)\cap S_{i,i+1}(A).$$
We have an inclusion $M^{A_i}\hk S_{i-1,i}(A)\cap S_{i,i+1}(A)$. (In fact,
according to the formula $(\ref{genvec})$, the cyclic vector $(\ref{vec})$
of $M^{A_i}$ belongs both to $S_{i-1,i}(A)$ and $S_{i,i+1}(A)$).
In addition, $\dim M^{A_i}=\dim (S_{i-1,i}(A)\cap S_{i,i+1}(A))$. To
prove the latter, recall that
$S_{i-1,i}(A)+ S_{i,i+1}(A)=S_{i-1,i+1}(A)$ (see lemma $(\ref{j-i})$). Thus
\begin{multline*}
\dim \left(S_{i-1,i}(A)\cap S_{i,i+1}(A)\right)=
\dim S_{i-1,i}(A)+\dim S_{i,i+1}(A)-\dim S_{i-1,i+1}(A)=\\
=\left(\prod_{j\ne i-1,i,i+1} a_j\right)\times
[ (a_i-a_{i-1}+1)a_{i+1}+a_{i-1}(a_{i+1}-a_i+1)-\\
-(a_{i+1}-a_{i-1}+1)a_i] =
\left(\prod_{j\ne i-1,i,i+1} a_j\right)\times (a_{i+1}-a_i+a_{i-1})=
\dim M^{A_i}.
\end{multline*}
Proposition is proved.
\end{proof}

Proposition $(\ref{filt})$ and remark $(\ref{stop})$ allows us to construct
our filtration. In fact, we
have a submodule $M^{A_i}$ in $S_{i,i+1}(A)$ and the quotient is isomorphic to
$S_{i,i+1}(A_{i-1,i})$. But this module also contains a corresponding fusion
as a submodule. So we can apply our lemma one more time and so on.
To formulate the condition of the finishing of our procedure, denote
by $s$ the map $\N^n\to\N^n$: $s(A)=A_{i-1,i}$. Let $s(A)_j, j=1,\ld,n$
be the elements of $s(A)$ with $s(A)_j\le s(A)_{j+1}$.
Our procedure will stop at the $k$-th step, if one of the below conditions
hold:
\begin{enumerate}
\item $s^k(A)_j=1$ for $j=1,\ld,i-1$. Then according to the first part of the
remark $(\ref{stop})$ $S_{i,i+1}(s^k(A))$ is some fusion product.
\item $s^k(A)_i=s^k(A)_{i+1}$. Then according to the second part of the remark
$(\ref{stop})$ $S_{i,i+1}(s^k(A))$ is also some fusion product.
\end{enumerate}
Let us give an example.\\
{\bf Example.} $A=(4,5,6,9), i=3$. We have:
\begin{gather*}
M^{(4,5,6,9)}/S_{3,4}(4,5,6,9)\simeq M^{(4,5,5,10)};\\
M^{(4,8)}\hk S_{3,4}(4,5,6,9),\qquad
S_{3,4}(4,5,6,9)/M^{(4,8)}\simeq S_{3,4}(4,4,7,9);\\
M^{(4,6)}\hk S_{3,4}(4,4,7,9),\qquad
S_{3,4}(4,4,7,9)/M^{(4,6)}\simeq S_{3,4}(3,4,8,9);\\
M^{(3,5)}\hk S_{3,4}(3,4,8,9),\qquad
S_{3,4}(3,4,8,9)/M^{(3,5)}\simeq S_{3,4}(3,3,9,9);\\
S_{3,4}(3,3,9,9)\simeq M^{(3,3)}.
\end{gather*}
Thus, we have a filtration of $M^{(4,5,6,9)}$ with the following quotients:\\
$M^{(4,8)}, M^{(4,6)}, M^{(3,5)}, M^{(3,3)}$ and $M^{(4,5,5,10)}$.

\subsection{Second description.}
For the second description we need an extra knowledge about
the structure of the vectors $w_j$ from $(\ref{genvec})$.
We will use a fermionic
realization of the fusion product from \cite{mi2}.

Let $F$ be the space of the semi-infinite forms (see \cite{mi2} or the
first section of this paper).
$F$ carries a structure of the $\slth$ module of the level $1$. There
exists a set of extremal vectors $v(i)\in F, i\in\Z$, such that
$e_{k-1} v(k)=v(k-2)$ and $e_{>k-1}v(k)=0$. Note also that
$U(\slth)\cdot v(0)\simeq L_{0,1}$ and $U(\slth)\cdot v(1)\simeq L_{1,1}$
(here $L_{0,1}$ and $L_{1,1}$ are two level $1$ irreducible $\slth$ modules).

In \cite{mi2} it was shown that we can embed the space $M^A$ into the
tensor power $F^{\T (a_n-1)}$. Let $d_i$ be the number of such $k$
that $a_k=i$. Suppose that $d_1=0$ (recall that
$M^{(1,a_2,\ld,a_n)}\simeq M^{(a_2,\ld,a_n)}$). Denote
$$v_A=v(n)\T v(n-d_2)\T v(n-d_2-d_3)\T \ld \T v(n-d_2-\ld -d_{a_n-1}).$$
Then $\Cn\cdot v_A\simeq M^A$. Another important fact is the following equality:
if $a_i\ne a_{i+1}$ then
\begin{multline}
\label{next}
[e_{(n)}(z)^{a_i-1}]_{N_A(a_i-1)}v_A \qquad\text{ is proportional to }\\
(e_{n-1}\T e_{n-d_2-1}\T\ld\T e_{n-d_2-\ld-d_{a_i-1}-1}\T \Id\T\ld\T \Id)
v_A=\\
=v(n-2)\T v(n-d_2-2) \T\ld\T v(n-d_2-\ld-d_{a_i-1}-2)\T\\
\T v(n-d_2-\ld-d_{a_i})\T \ld\T v(n-d_2-\ld-d_{a_n-1}).
\end{multline}
We generalize the formula $(\ref{next})$ in the following way: let $k$ be
such number that $a_i-1< k< a_{i+1}-1$. Note that
\begin{multline*}
v_A=\left[v(n)\T v(n-d_2)\T\ld\T v(n-d_2-\ld-d_{a_i-1})\right]\T\\
\T\left[v(n-d_2-\ld-d_{a_i})^{\T (a_{i+1}-a_i)}\right] \T\\
\T\left[v(n-d_2-\ld-d_{a_{i+1}})\T\ld\T v(n-d_2-\ld-d_{a_n-1})\right].
\end{multline*}
Let
\begin{gather}
\label{v123}
v_1=v(n)\T v(n-d_2)\T\ld\T v(n-d_2-\ld-d_{a_i-1});\\
v_2=v(n-d_2-\ld-d_{a_i})^{\T a_{i+1}-a_i};\nonumber\\
v_3=v(n-d_2-\ld-d_{a_{i+1}})\T\ld\T v(n-d_2-\ld-d_{a_n-1}).\nonumber
\end{gather}
\begin{lem}
Introduce a notation
\begin{gather}
\label{u1}
u_1=(e_{n-1}\T e_{n-d_2-1}\T\ld\T e_{n-d_2-\ld-d_{a_i-1}-1})v_1,\\
\label{u2}
u_2(k)=(e_{n-d_2-\ld-d_{a_i}-1}^{k-a_i+1})v_2,
\end{gather}
where $e_{n-d_2-\ld-d_{a_i}-1}^{k-a_i+1}$ is a power of the
operator $e_{n-d_2-\ld-d_{a_i}-1}$, which acts on the space
$F^{\T (a_{i+1}-a_i)}$ as
on the tensor power of the Lie algebra module.
Then
$[e_{(n)}(z)^k]_{N_A(k)}v_A$ is proportional to the vector
\begin{equation}
\label{u12}
u_1 \T u_2(k)\T v_3.
\end{equation}
\end{lem}
\begin{proof}
One can prove this lemma in the same way, as it was done in \cite{mi2}
in the case $k=a_i-1$.
\end{proof}
Now we can describe the structure of the module $S_{i,i+1}(A)$. Recall that
\begin{equation}
\label{gen}
S_{i,i+1}(A)=\Cn\cdot \bra w_{a_i-1},\ld , w_{a_{i+1}-1} \ket
\end{equation}
and
$w_j=[e_{(n)}(z)^j]_{N_A(j)}v_A.$
From the above lemma we obtain
\begin{multline*}
w_j=\const \cdot v(n-2)\T\ld\T v(n-d_2-\ld-d_{a_i-1}-2)\T\\
\T [(e_{n-d_2-\ld-d_{a_i}-1}\T\underbrace{\Id\T\ld\T\Id}_{a_{i+1}-a_i-1}+\ld
+\underbrace{\Id\T\ld\T\Id}_{a_{i+1}-a_i-1}\T
e_{n-d_2-\ld-d_{a_i}-1})^{j-a_i+1} v_2]\T\\\T v_3.
\end{multline*}
Hence we can rewrite the formula $(\ref{gen})$ in the following way:
\begin{prop}
\label{mprop}
Denote
$$A'=(a_1,\ld,a_{i-1},\underbrace{a_i,\ld,a_i}_{n-i-1});\qquad
A''=(a_{i+1}-a_i+1, a_{i+2}-a_i+1,\ld, a_n-a_i+1).$$
Then
$S_{i,i+1}(A)\hk M^{A'}\T M^{A''}$ and the image of $S_{i,i+1}(A)$ is
generated from the vector
$v_{A'}\T v_{A''}$ by the action of
$\C[e_0,\ld,e_{n-3}, e^{(2)}_{n-i-1}]$, where the operator
$e^{(2)}_{n-i-1}$ acts on $M^{A'}\T M^{A''}$ as
$\Id\T e_{n-i-1}$.
\end{prop}
\begin{proof}
Note that for $u_1, u_2(k)$ defined by $(\ref{u1}), (\ref{u2})$ and $v_3$
from the formula $(\ref{v123})$ we have
$u_1\in M^{A'}$ and $u_2(k)\T v_3\in M^{A''}$. Moreover,
$w_{a_i-1}=v_1$ and $w_j=(e^{(2)}_{n-i-1})^{j-a_i+1} w_{a_i-1}$.
Thus $S_{i,i+1}(A)$ is a submodule
of the module, described in the proposition. To prove that we have an
equality, it is enough to mention that
$(e^{(2)}_{n-i-1})^{a_{i+1}-a_i+1}=~0$.
\end{proof}

\begin{rem}
\label{alg}
Note that $e_{>n-i-1}$ vanishes on $M^{A''}$. Thus $S_{i,i+1}(A)$ is a
cyclic submodule of $M^{A'}\T M^{A''}$ with respect to the algebra,
generated by the operators
$$e^{(1)}_0+e^{(2)}_0,\ld, e^{(1)}_{n-i-2}+e^{(2)}_{n-i-2},
e^{(1)}_{n-i-1}, e^{(2)}_{n-i-1}, e^{(1)}_{n-i},\ld, e^{(1)}_{n-3}$$
(here $e_j^{(1)}=e_j\T\Id, e_j^{(2)}=\Id\T e_j$).
\end{rem}
\begin{rem}
\label{remtens}
Note that in the case $i=n-1$ we have
\begin{equation}
\label{tens}
S_{n-1,n}(A)\simeq M^{(a_1,\ld,a_{n-2})}\T \C^{a_n-a_{n-1}+1}.
\end{equation}
In fact, proposition $(\ref{mprop})$ gives us an embedding of the left hand
side of $(\ref{tens})$ into the right hand side. But we have an operator
$e^{(2)}_{0}$ acting on $S_{n-1,n}(A)$. Thus
\begin{multline*}
S_{n-1,n}(A)\hookleftarrow \C[e^{(1)}_{n-3},\ld,e^{(1)}_0,e^{(2)}_0]\cdot
(v_{(a_1,\ld,a_{n-2})}\T v_{a_n-a_{n-1}+1}) =\\
M^{(a_1,\ld,a_{n-2})}\T \C^{a_n-a_{n-1}+1}.
\end{multline*}
\end{rem}

Thus we have an embedding of $S_{i,i+1}(A)$ into the tensor product
$M^{A'}\T M^{A''}$.
Recall that $A'$ contains a lot of
coinciding
elements. By some reasons, discussed in the previous section, in this
work we concentrate on those sets, which have no coinciding elements.
Thus we will construct another embedding of $S_{i,i+1}(A)$, which looks
somehow strange, but its usefulness will be explained in the sequel
section.
\begin{prop}
\label{emb}
Let $A=(a_1<\ld <a_n)$. Suppose that for any $j>i$ we have
$a_{j+1}-a_j>1$.
Introduce a notation:
\begin{gather*}
A_1=(a_1,\ld,a_{i-1}, a_i+1, a_i+2,\ld, a_i+n-i-1);\\
A_2=(a_{i+1}-a_i+1, a_{i+2}-a_i, a_{i+3}-a_i-1, \ld, a_n-a_i-n+i+2).
\end{gather*}
(Note that the elements of $A_1$ and $A_2$ increase). Then we have an
embedding
$S_{i,i+1}(A)\hk M^{A_1}\T M^{A_2}$ and $S_{i,i+1}(A)$ is generated from
$v_{A_1}\T v_{A_2}$ by the action of the algebra
$\C[e_0,\ld,e_{n-3}, e^{(2)}_{n-i-1}]$.
\end{prop}
\begin{proof}
In our case
\begin{gather*}
v_1=v(n)^{\T (a_1-1)} \T v(n-1)^{\T (a_2-a_1)}\T\ld\T
v(n-i+1)^{\T (a_i-a_{i-1})};\\
u_1=v(n-2)^{\T (a_1-1)} \T v(n-3)^{\T (a_2-a_1)}\T\ld\T
v(n-i-1)^{\T (a_i-a_{i-1})};\\
v_2=v(n-i)^{\T (a_{i+1}-a_i)};\\
v_3=v(n-i-1)^{\T (a_{i+2}-a_{i+1})}\T \ld\T v(1)^{\T (a_n-a_{n-1})}.
\end{gather*}
We can rearrange the factors in the tensor product $u_1\T v_2\T v_3$ to
obtain the product $u_1'\T v_2\T v_3'$, where
\begin{multline*}
u_1'= v(n-2)^{\T (a_1-1)} \T \ld\T
v(n-i-1)^{\T (a_i-a_{i-1})}\T\\ \T v(n-i-1)\T v(n-i-2)\T\ld\T v(1);\\
\shoveleft
{v_3'=v(n-i-1)^{\T (a_{i+2}-a_{i+1}-1)}\T v(n-i-2)^{\T (a_{i+3}-a_{i+2}-1)}\ld
\T v(1)^{\T (a_n-a_{n-1}-1)}.}
\end{multline*}
Note that
$$M^{A_1}\simeq \C[e_0,\ld,e_{n-3}]\cdot u_1',\quad
M^{A_2}=\C[e_0,\ld, e_{n-i-1}]\cdot (v_2\T u_3').$$
To finish the proof of the proposition it is enough, to recall that
$S_{i,i+1}(A)$ is generated from the vectors $w_j$ and
$w_j=u_1'\T (e_{n-i-1}^{j-a_i+1} v_2)\T v_3'.$
\end{proof}

We finish this subsection with the discussion of the Lie algebra, acting on
$S_{i,i+1}(A)$.
Let $\Lie_{i,n}=\slt\T B_{i,n}$, $0< i\le n$, where
$B_{i,n}$ is a
commutative graded associative algebra with a generators $t$ of degree $1$
and $u$ of degree $n-i$ and relations $t^n=0, u^2=0, tu=0$. For $i=0$ let
$\Lie_{0,n}=\slt\T (\C[t]/t^{n+1})$.
From the propositions $(\ref{mprop}), (\ref{emb})$ we obtain the following
lemma:
\begin{lem}
\label{cycvec}
$S_{i,i+1}(A)$ is cyclic $\Lie_{i-1,n-2}$ module with a cyclic vector being
the tensor product
$v_{A'}\T v_{A''}$ (see proposition $(\ref{mprop}$)) or
$v_{A_1}\T v_{A_2}$ (see proposition $(\ref{emb}$)).
We denote this cyclic vector by $v_{i,i+1}(A)$.
\end{lem}
For example, for $i=1$ we obtain that $S_{1,2}(A)$ is
$\slt\T (\C[t]/t^{n-1})$ module.
But as it was mentioned above, $S_{1,2}(A)\simeq M^{(a_2-a_1+1,a_3,\ld,a_n)}$.

Note that $S_{n-1,n}(A)$, being a tensor product
$M^{(a_1,\ld,a_{n-2})}\T \C^{a_n-a_{n-1}+1}$, is a cyclic
$\slt\T (\C[t]/t^{n-2})$ module, but
its cyclic vector is not the tensor product of the lowest weight vectors of
$M^{(a_1,\ld,a_{n-2})}$ and $\C^{a_n-a_{n-1}+1}$.  To be specific,
one can put this cyclic vector $l$ to be
$u_{(a_1,\ld,a_{n-2})}\T v_{a_n-a_{n-1}+1}$ (recall that $u_A$ is a highest
vector in $M^A$ with respect to the $h_0$ grading).
\begin{lem}
\label{e1}
$l=e_1^{a_1+\ld +a_{n-1}-n} v_A.$
\end{lem}
\begin{proof}
This follows from the fermionic realization of $M^A$. Note also that
$e_1 l=~0$.
\end{proof}

\subsection{Third description: induction.}
In this section we obtain the analogue of the lemma $(\ref{demazure})$
for the modules $S_{i,i+1}(A)$. In what follows we consider the module
$M^{(a_1,\ld,a_{n-1})}$ as a submodule of $M^A$ via the isomorphism
$$M^{(a_1,\ld,a_{n-1})}\simeq \C[e_1,\ld,e_{n-1}]\cdot v_A.$$
\begin{lem}
\label{ind1}
Let $i\ne n-1$. Then we have (see lemma $(\ref{cycvec})$):
$$v_{i,i+1}(A)=v_{i,i+1}(a_1,\ld,a_{n-1}).$$
It means that
$$S_{i,i+1}(a_1,\ld,a_{n-1})\simeq
\C[e\T t,\ld,e\T t^{n-3}, e\T u]\cdot v_{i,i+1}(A)$$
and
$$ S_{i,i+1}(A)=\C[e_0]\cdot S_{i,i+1}(a_1,\ld,a_{n-1}).$$
\end{lem}
\begin{proof}
This is an immediate consequence from the construction of the embedding from
the proposition $(\ref{mprop})$.
\end{proof}

Now let $i=n-1$. Consider the operator $e_1$ acting
on $M^A$. Note that $e_1^{a_1+\ld +a_{n-1}-n+1}=0$. Let
$l=e_1^{a_1+\ld +a_{n-1}-n}v_A$. From one hand, $l$ is a highest weight
vector from $M^{(a_1,\ld,a_{n-1})}$. From the other hand, $l$ is a cyclic
vector in $S_{n-1,n}(A)$ with respect to the algebra $\C[e_0,\ld,e_{n-3}]$
(see  lemma $(\ref{e1})$). Thus we obtain the following lemma:
\begin{lem}
\label{ind2}
There is a vector $l$ ,
$$l\in S_{n-1,n}(A)\simeq M^{(a_1,\ld, a_{n-2})}\T \C^{a_n-a_{n-1}+1}\hk M^A$$
such that $l$ is a cyclic vector of $S_{n-1,n}(A)$ with respect to the algebra
$\C[e_0,\ld,e_{n-3}]$ and $l=u_{(a_1,\ld,a_{n-1})}$.
In particular it means that
$$
S_{n-1,n}(A)=\C[e_0]\cdot (\C[f_0,\ld, f_{n-2}]\cdot u_{(a_1,\ld,a_{n-1})}).
$$
\end{lem}

\section{Algebro-geometric properties of the Schubert varieties.}
The main point of this section is the study of the line bundles on $\shn$.
We start with the description of $\shn$ as a projective algebraic variety.
\subsection{$\shn$ as an algebraic variety.}
Recall (see \cite{mi1}) that the dual space to the module $M^A$,
$A=(a_1\le\ld\le a_n)$ can be realized as the space
of the symmetric polynomials $f(z_1,\ld, z_s), s\ge 0$ with the following
conditions: $\deg_{z_i}f<n$ and
$$f(\underbrace{z,\ld,z}_i,z_{i+1},\ld,z_s)\div
z^{\sum_{j=1}^n (i+1-a_j)_+},\quad i=1,2,\ld,$$
where $a_+=\max (a,0)$ and for the polynomials $p$ and $q$ $p\div q$ means
that $q$ divides $p$.
(Note that $(M^A)^*$ is naturally graded by the action of $h_0$ and thus
$(M^A)^*=\bigoplus_s (M^A)^*(s)$. Then $(M^A)^*(s)$ is realized in the space
of the polynomials in $s$ variables with the above condition).
\begin{lem}
\label{dual}
Let $A=(a_1\le\ld\le a_n)\in\N^n, B=(b_1\le\ld\le b_n)\in\N^n$ and
$C=(a_1+b_1-1,\ld, a_n+b_n-1)$. Then the following is true:\\
$1).$\ There is an embedding of $\C[e_0,\ld,e_{n-1}]$ modules
$\alpha: M^C\hk M^A\T M^B$  sending $v_C$ to
$v_A\T v_B$ (in fact, $\al$ is $\sltc$ homomorphism).\\
$2).$\ There is a surjection
$\be: (M^A)^*\T (M^B)^*\sur (M^C)^*$ given by the formula
$$f(z_1,\ld,z_{s_1})\T g(z_1,\ld,z_{s_2})\mapsto h(z_1,\ld,z_{s_1+s_2}),$$
where
$$h(z_1,\ld,z_{s_1+s_2})=\sum f(z_{\sigma(1)},\ld,z_{\sigma({s_1})})
g(z_{\tau(1)},\ld,z_{\tau({s_2})}),$$
and the sum is taken over such pairs $(\sigma,\tau)$
$$\sigma: \{1,\ld,s_1\}\to \{1,\ld,s_1+s_2\},\
\tau: \{1,\ld,s_2\}\to \{1,\ld,s_1+s_2\},$$ that
$\sigma(i)<\sigma (i+1), \tau (j)<\tau (j+1)$ and the images of $\sigma$
and $\tau$ do not intersect.\\
$3).$\ $\al^*=\be$.
\end{lem}
\begin{proof}
First statement is a special case of the proposition $(\ref{tg})$. One can
check
that the image of $\be$ is really a subspace of $(M^C)^*$. To
prove that $\beta$ is a surjection it is enough to show $3).$  To make it
clear, recall that the isomorphism between the dual space $(M^A)^*$ and the
above space of the symmetric polynomials has the following form:
$$(M^A)^*\ni \theta\mapsto
f_{\theta}=\sum_{i_1,\ld,i_s} z_1^{i_1}\ld z_s^{i_s}
\theta (e_{n-1-i_1}\ld e_{n-1-i_s}v_A).$$
This finishes the proof.
\end{proof}

Introduce a notation:
$$A(k)=(ka_1-k+1, ka_2-k+1,\ld, ka_n-k+1),\quad k=1,2,\ld.$$
For example, $A(1)=A$. From the above lemma we obtain a graded algebra
$F_A=\bigoplus_{i=0}^{\infty} (M^{A(i)})^*$ (we put $(M^{A(0)})^*=\C$).
Note that
$F_A$ is generated by its first degree component $(M^{A(1)})^*=(M^A)^*$.

Now fix some basis $v_1,\ld,v_N$ in $M^A$ (surely,
$N=\prod_{i=1}^n a_i$). Denote also by $\xi_i$ the dual basis in $(M^A)^*$.
Recall that $\shn\hk \Pro(M^A)$. Thus, any point ${\bf x}\in\shn$ can be
written in the basis $v_i$: ${\bf x}=(x_1:\ld :x_N)$.
\begin{prop}
\label{coordring}
$F_A$ is the coordinate ring of $\shn$.
\end{prop}
\begin{proof}
We need to prove that for any homogeneous polynomial $p$ in $N$ variables
the following
two statements are equivalent:
\begin{enumerate}
\item for any ${\bf x}\in\shn$ $p(x_1,\ld,x_N)=0$.
\item $p(\xi_1,\ld,\xi_N)=0$ (here $\xi_i$ and $p(\xi_1,\ld,\xi_N)$ are
considered as an elements of the algebra $F_A$).
\end{enumerate}

First, let $p(\xi_1,\ld,\xi_N)=0$. Then for any $w\in M^{A(N)}$
$p(\xi_1,\ld,\xi_N)w=0$. From the lemma $(\ref{dual})$ we obtain an
embedding: $\varphi: M^{A(N)}\to (M^A)^{\T N}$. Let
$w=\exp(\sum_{i=0}^{n-1} e_it_i)\cdot v_{A(N)}$, where $t_i$ are some complex
numbers. Note that $\varphi v_{A(N)}=v_A^{\T N}$ and the operator $e_k$,
acting on $M^{A(N)}$ is a sum $\sum_{i=1}^N e_k^{(i)}$ of the operators,
acting on the corresponding copies of $M^{A}$. Thus
\begin{multline*}
w=\exp\left(\sum_{i=0}^{n-1} e_it_i\right)\cdot v_{A_N}= \\
=\exp\left(\sum_{i=0}^{n-1} \left(\sum_{j=1}^N e_i^{(j)}\right) t_i\right)
\cdot v_A^{\T N}=
\bigotimes_{j=1}^N \left[\exp\left(\sum_{i=0}^{n-1} e_i^{(j)} t_i\right)\cdot
v_A\right].
\end{multline*}
Now, let $\exp\left(\sum_{i=0}^{n-1} e_i^{(j)} t_i\right)\cdot v_A=
\sum_{i=1}^N x_iv_i$. Then
$$w=\bigotimes_{j=1}^N \left( \sum_{i=1}^N x_iv_i \right)$$
and because of the part $3)$ of the lemma
$(\ref{dual})$ we obtain
$$0=p(\xi_1,\ld,\xi_N)w=p(x_1,\ld,x_N).$$
Hence we have proved that $p({\bf x})=0$ for ${\bf x}$ from the
big cell $U_x=\{\exp \left(\sum t_ie_i\right)\cdot [v_A]\}$. But $\shn$ is a
closure of the big cell. Thus $p$ vanishes on all $\shn$.

Now, let $p(x_1,\ld,x_N)=0$ for any $(x_1:\ld :x_N)\in\shn$. Then $p$
vanishes on the big cell. The same considerations as in the first part of the
proof give us $p(\xi_1,\ld,\xi_N)=0$.
\end{proof}

\subsection{Line bundles on $\shn$.}
It was proved in the first section that there is a chain of bundles
\begin{equation}
\label{fibr}
\shn\to \sh^{(n-1)}\to\ld\to \sh{(1)}\simeq\Pro^1.
\end{equation}
Thus $\Pic(\shn)=\Z^n$. We want to fix this isomorphism. Let
$C_i, i=0,\ld,n-1$ be a collection of projective lines in $\shn$. Namely,
$$C_i=\ov{\left\{\exp (e_it)\cdot [v_A], t\in\C\right\}}$$
(we fix some $A$ with no coinciding
elements and identify $\shn$ with $\sh_A$. Surely, our lines do not depend
on the choice of $A$). For example, $C_0=\Slt\cdot [v_A]$. Because of the
fibrations $(\ref{fibr})$, any line bundle is determined by its
restriction on the lines $C_i$. Let $\E$ be a line bundle, such that
$\E|_{C_i}\simeq \Ob(b_i)$. Then we denote this bundle as
$$\E=\Ob(b_{n-1},b_{n-2}-b_{n-1},b_{n-3}-b_{n-2},\ld,b_0-b_1).$$

Recall that for any $A$ with $a_i\ne a_j$ we have an embedding
$\imath_A: \shn\hk \Pro(M^A)$. Thus for any such $A$ there is a bundle
$\imath_A^* \Ob(1)$ on  $\shn$. The following lemma explains the
choosed parameterization of the set of bundles.
\begin{lem}
$\imath_A^* \Ob(1)\simeq \Ob(a_1-1,a_2-1,\ld,a_n-1)$.
\end{lem}
\begin{proof}
Note that $H^0(\shn,\E|_{C_i})\simeq \left(\C[e_i]\cdot v_A\right)^*$.
Recall that for any $A$ we have an isomorphism (see lemma
$(\ref{demazure})$):
\begin{equation}
\C[e_i,\ld,e_{n-1}]\cdot v_A\simeq M^{(a_1,\ld,a_{n-i})},
\quad e_j\mapsto e_{j-i}.
\end{equation}
Thus $\C[e_i]\cdot v_A$ is $(a_1+\ld+a_{n-i}-n+i+1)$ -- dimensional space,
and so $(\imath^* \Ob(1))|_{C_i}\simeq\Ob(a_1+\ld+a_{n-i}-n+i)$.
Lemma is proved.
\end{proof}

Recall that there is a bundle $\pi_{n,n-1}:\shn\to \sh^{(n-1)}$
with a fiber $\Pro^1$. The following lemma shows that there
exists a rank two bundle $\xi$ on $\sh^{(n-1)}$ such that $\shn$ is a
projectivization of $\xi$.
\begin{lem}
$\shn\simeq \Pro\left((\pi_{n,n-1})_* \imath_{(2,3,\ld,n+1)}^* \Ob(1)\right)$.
\end{lem}
\begin{proof}
It follows from the fact that
$\imath_{(2,3,\ld,n+1)}^* \Ob(1)|_{\pi_{n,n-1}^{-1} x}=\Ob(1)$ for any
$x\in \sh^{(n-1)}$.
\end{proof}
{\bf Example.} In the case $n=2$ we have
$\sh^{(2)}\simeq \Pro(\Ob(0)\oplus \Ob(2))$.

We finish this subsection with the computation of the cannonical bundle
$K_{\shn}$.
\begin{lem}
$K_{\shn}\simeq \Ob(\underbrace{-2,\ld,-2}_n)$.
\end{lem}
\begin{proof}
We prove lemma by the induction on $n$. For $n=1$ $\sh^{(1)}\simeq\Pro^1$ and
$K_{\Pro^1}\simeq\Ob(-2)$. Suppose our lemma is true for $n-1$. Then because
of the existence of the bundle $\pi_{n,1}:\shn\to\Pro^1=C_0$ we obtain that
$K_{\shn}=\Ob(\underbrace{-2,\ld,-2}_{n-1},a)$ for some $a\in\Z$. Recall
that there is an
open (in the Zariski topology) set $\G\cdot [v_A]\hk\shn$ ($\G$ is a group
$\Slt (\C[t]/t^n)$),
which is fibered over $C_0$ with a fiber $\C^{n-1}$. Recall also
(see lemma $(\ref{bunlem})$) that we have an equality
$$\exp\left(\sum_{i=0}^{n-1} x_ie_i\right)\cdot [v_A]=
  \exp\left(\sum_{i=0}^{n-1} y_if_i\right)\cdot [u_A]$$
if and only if for two polynomials $x(t)=\sum_{i=0}^{n-1} x_it^i$ and
$y(t)=\sum_{i=0}^{n-1} y_it^i$ we have $x(t)y(t)=1$ in the ring
$\C[t]/t^n$. To compute the restriction of $K_{\shn}$ to $C_0$ we need
to rewrite the form $dy_0\wedge dy_1\wedge \ld\wedge dy_{n-1}$ in the
$x_i$-coordinates. It is easy to check that the result will be the following
$n$-form: $\frac{(-1)^n}{x_0^{2n}} dx_0\wedge \ld\wedge dx_{n-1}$. Thus the
restriction of $K_{\shn}$ to $C_0$ is $\Ob(-2n)$. Using the fact that
$K_{\shn}=\Ob(-2,\ld,-2,a)$ we obtain $a=-2$. Lemma is proved.
\end{proof}

\subsection{Fusion products as a dual spaces of sections of the line bundles.}
The main goal of this subsection is to prove the following theorem:
let $A=(a_1,\ld,a_n), 1\le a_1\le\ld\le a_n$. Then the dual space
$H^0 (\Ob(a_1-1,\ld,a_n-1))^*$ as $\sltc$ module is isomorphic to the
fusion product $M^A$.

Let $\G=\Slt(\C[t]/t^n)$.
Recall the submodule $S_{i,i+1}(A)\hk M^A$ which is the cyclic
$\Lie_{i-1,n-2}$ module  with a cyclic vector $v_{i,i+1}(A)$. Denote by
$\Lg_{i,n}$ the Lie group of the Lie algebra $\Lie_{i,n}$.
\begin{opr}
The variety $\sh^{(n-1)}_i\hk \Pro (S_{i,i+1}(A))\hk \Pro (M^A)$ is a closure
of the orbit of the
point $[v_{i,i+1}(A)]:$
$$\sh_i^{(n-1)}=\ov{\Lg_{i-1,n-2}\cdot [v_{i,i+1}(A)]}.$$
\end{opr}
\begin{lem}
$1).$\ $\sh_i^{(n-1)}$ is a subvariety of $\shn$.\\
$2)$.\ $\shn=\G\cdot [v_A]\cup \bigcup_{i=1}^{n-1} \sh_i^{(n-1)}$.
\end{lem}
\begin{proof}
We will prove both statements together by the induction on $n$. Consider
a fiber $\pi_{n,1}^{-1}\pi_{n,1}([v_A])\simeq\sh^{(n-1)}$ of the bundle
$\pi_{n,1}:\shn\to\Pro^1$.
By the induction assumption
we have the decomposition of this fiber into a cell
$U'_x=\left\{\exp\left(\sum_{i=1}^{n-1} e_it_i\right)\cdot [v_A]\right\}$,
divisors
$\sh_i^{(n-2)}, i=1,\ld,n-2$ and
$$M=\ov{\left\{\pi_{n,1}^{-1}\pi_{n,1}([v_A])\setminus U'_x
\setminus\bigcup_{i=1}^{n-2} \sh_i^{(n-2)}\right\}}\simeq \sh^{(n-2)}.$$
First note that $\G\cdot [v_A]=\Slt\cdot U'_x$.
Next, because of the lemma $(\ref{ind1})$ we have
$$\sh_i^{(n-1)}=\Slt\cdot \sh_i^{(n-2)},\qquad i=1,\ld,n-2.$$
Finally, because of the lemma $(\ref{ind2})$ we obtain
$\sh_{n-1}^{(n-1)}=\Slt\cdot M$. Lemma is proved.
\end{proof}
The following proposition gives a description of the varieties
$\sh^{(n-1)}_i$.
For the convenience, let $\sh^{(0)}$ be a point and
$\pi_{n,0}$ a unique map $\shn\to \sh^{(0)}$.
\begin{prop}
\label{descr}
Let $i=2,\ld,n-1$. Then
$$\sh_i^{(n-1)}=\{(x,y)\in \sh^{(n-2)}\times \sh^{(n-i)}:\
\pi_{n-2,n-i-1}(x)=\pi_{n-i,n-i-1}(y)\}$$
(recall that $\pi_{n,k}$ is a projection $\shn\to\sh^{(k)}$ with a fiber
$\sh^{(n-k)}$). In the case $i=1$ we have $\sh_1^{(n-1)}\simeq \sh^{(n-1)}$.
\end{prop}
\begin{proof}
Note that $S_{1,2}(A)\simeq M^{(a_2-a_1+1,a_3,\ld,a_n)}$. Thus
$\sh_1^{(n-1)}\simeq \sh^{(n-1)}$.

Recall that we have proved in the proposition $(\ref{emb})$ that for any
$i=2,\ld,n-1$ and some special $B\in\N^{n-2}, C\in\N^{n-i}$ with no coinciding
elements in both $B$ and $C$ we have an embedding
$$S_{i,i+1}(A)\hk M^B\T M^C,\quad v_{i,i+1}(A)\mapsto v_B\T v_C.$$
In addition, $S_{i,i+1}(A)$ is isomorphic to
\begin{multline}
\label{formula}
\C[e^{(1)}_0+e^{(2)}_0,\ld, e^{(1)}_{n-i-2}+e^{(2)}_{n-i-2},
e^{(1)}_{n-i-1}+e^{(2)}_{n-i-1},\\ e^{(2)}_{n-i-1}, e^{(1)}_{n-i}, \ld,
e^{(1)}_{n-3}]\cdot (v_B\T v_C),
\end{multline}
where the upper index $(i)$ means that the corresponding operator is acting
on the $i$-th factor of the tensor product $M^B\T M^C$.
This shows that $\sh^{(n-1)}_i$ is embedded to the Cartesian product
$\sh^{n-2}\times \sh^{n-i}$. Moreover, formula $(\ref{formula})$ shows that
points from $\Ga_{i,n}\cdot [v_{i,i+1}(A)]$ (where
$\Ga_{i,n}$ is the Lie group of an abelian Lie algebra, spanned by
$e^{(1)}_0+e^{(2)}_0,\ld,$ $e^{(1)}_{n-i-2}+e^{(2)}_{n-i-2},
e^{(1)}_{n-i-1}+e^{(2)}_{n-i-1},$ $e^{(2)}_{n-i-1}, e^{(1)}_{n-i}, \ld,
e^{(1)}_{n-3})$
are the pairs
$(x,y)$ from the product of the big cells of $\sh^{(n-2)}$ and $\sh^{(n-i)}$
with a property
$\pi_{n-2,n-i-1}(x)=\pi_{n-i,n-i-1}(y)$. The equality
$\sh^{(n-1)}_i=\ov{\Ga_{i,n}\cdot [v_{i,i+1}(A)]}$ finishes the proof.
\end{proof}
\begin{cor}
In the case $i=n-1$ we have $\sh^{(n-1)}_{n-1}\simeq\Pro^1\times \sh^{(n-2)}$
(see the remark $(\ref{remtens}))$.
\end{cor}

Our goal is to prove that for $A=(a_1\le\ld\le a_n)$ fusion product
$M^A$ is realized in the dual space of the sections of $\Ob(a_1-1,\ld,a_n-1)$.
The way of the proof is the restriction of
the bundles to the subvarieties $\sh_i^{(n)}$. We need some additional
information about $\sh_i^{(n-1)}$.

Recall that $\shn$ is fibered over $\Pro^1$ with a fiber $\sh^{(n-1)}$.
By the induction on $n$ one can show that $\shn$ is a union of
the cells $\C^i$ (which do not intersect) and the number
of the cells $\C^i$
is equal to $\binom{n}{i}$. Thus we obtain
$H_{2i}(\shn,\Z)=\Z^{l_i}$, $l_i=\binom{n}{i}$ and $H_{2i-1} (\shn,\Z)=0$.
Recall
that we have defined the set of projective lines $C_i$
($C_i=\ov{\{\exp(te_i) \cdot [v_A]\}}$).
One can show that the classes of $C_i$ in
$H_2(\shn,\Z)$ are generators of the latter group. For the variety
$M\hk\shn$ of the complex dimension $m$ we write $[M]$ for the corresponding
class in $H_{2m}(\shn,\Z)$.

Denote by $\sh_n^{(n-1)}$ the fiber $\pi_{n,1}^{-1}\pi_{n,1}([u_A])$:
$\sh_n^{(n-1)}=
\ov{\left\{\exp\left(\sum_{i=1}^{n-1} t_i f_i\right)\cdot [u_A]\right\}}.$
Surely, $\sh_n^{(n-1)}\simeq \sh^{(n-1)}$.
Thus, we have a varieties $\sh_j^{(n-1)}$ for $j=1,\ld,n$.
\begin{lem}
\label{intind}
$[C_i]\cdot [\sh_j^{(n-1)}]=\delta_{i,n-j}$, where $[M]\cdot [N]$ is an
intersection number of $M$ and $N$.
\end{lem}
\begin{proof}
First, we consider a special cases of $i=0,1$ and $j=n,n-1$. All other
cases are the consequences of these special ones.

Let $i=0$. Then $C_0$ doesn't
intersect with $\sh^{(n-1)}_i$ for $i<n$ and
$C_0\cap \sh_n^{(n-1)}=[u_A]$.

Let $j=n$. Then $[\sh_n^{(n-1)}]\cdot [C_i]=0$ for $i>0$ because
$\sh_n^{(n-1)}= \pi_{n,1}^{-1}\pi_{n,1} ([u_A])$ and
$C_i\hk \pi_{n,1}^{-1}\pi_{n,1} ([v_A])$ for $i>0$.

Let $i=1$. Then because of the lemma $(\ref{ind2})$ we have
$[C_1]\cdot [\sh_{n-1}^{(n-1)}]=1$. Note also that
$[C_1] \cdot [\sh_i^{(n-1)}]=0$ for $i<n-1$. In fact, because of the lemma
$(\ref{ind1})$ we have
$\sh_i^{(n-1)}\cap \pi_{n,1}^{-1}\pi_{n,1}([v_A])=\sh_i^{(n-2)}$. But $C_1$,
which
is embedded to the fiber $\pi_{n,1}^{-1}\pi_{n,1}([v_A])\simeq \sh^{(n-1)}$,
plays there a role of $C_0\hk\shn$. Thus, we are in the situation of the case
$i=0$.

Let $j=n-1$. We need to show that $\sh^{(n-1)}_{n-1}$ doesn't intersect
with $C_i, i>1$. Recall that
$\sh_{n-1}^{(n-1)}\cap \pi_{n,1}^{-1}\pi_{n,1}([v_A])=\sh_{n-1}^{(n-2)}$.
Thus we are in the situation of the case $j=n$.

All other cases can be considered in the same way, using lemmas
$(\ref{ind1},\ref{ind2})$.
\end{proof}
\begin{cor}
The classes $[\sh^{(n-1)}_i]$, $i=1,\ld,n$ are generators of the group
$H_{2n-2}(\shn,\Z)$. Moreover, via the identification
$H_{2n-2}(\shn,\Z)\simeq H^2(\shn,\Z)$
$\shin\mapsto [C_{n-i}]^*$.
\end{cor}

Now we need some fact about the restriction of line bundles on $\shn$ to
$\shin$.
Recall (see proposition $(\ref{descr})$) that there exists a bundle
$\nu_{n,i}:\shin\to \sh^{(n-i-1)}$ with a fiber $\Pro^1\times \sh^{(i-1)}$.
Thus, we can first restrict any bundle $\E$ on $\shn$ to
$\Pro^1\times \sh^{(i-1)}$ and then to $\Pro^1\times x,\ x\in \sh^{(i-1)}$.
Denote this restriction by $r_i(\E)$ (surely, $r_i(\E)$ doesn't depend on
$x$).
\begin{lem}
Let $\E=\Ob(a_1,\ld,a_n)$. Then $r_i(\E)=\Ob(a_{i+1}-a_i)$.
\end{lem}
\begin{proof}
Note that if the statement of the lemma is true for $\E_1$ and $\E_2$, then
it is also true for $\E_1^*$ and for $\E_1\T \E_2$. Hence it is enough to
prove our lemma in the case $0\le a_1<\ld <a_n$.

If $0\le a_1<\ld <a_n$, then $\E=\imath_A^*\Ob(1)$. Recall (see proposition
$(\ref{mprop})$) that $S_{i,i+1}(A)\hk M^{A'}\T M^{A''}$ and
\begin{multline*}
S_{i,i+1}(A)=\\=\C[e^{(1)}_0+e^{(2)}_0,\ld,
e^{(1)}_{n-i-1}+e^{(2)}_{n-i-1}, e^{(2)}_{n-i-1}, e^{(1)}_{n-i}, \ld,
e^{(1)}_{n-3}]\cdot (v_{A'}\T v_{A''}).
\end{multline*}
Thus we obtain
\begin{equation*}
H^0(r_i(\E))=(\C[e^{(2)}_{n-i-1}]\cdot (v_{A'}\T v_{A''}))^*.
\end{equation*}
But $(e^{(2)}_{n-i-1})^{a_{i+1}-a_i+1} (v_{A'}\T v_{A''})=0$ and
$(e^{(2)}_{n-i-1})^{a_{i+1}-a_i} (v_{A'}\T v_{A''})\ne 0$. Lemma is proved.
\end{proof}

\begin{cor}
\label{coheq}
Fix a number $i$: $1\le i < n$. Then for any $A=(a_1,\ld,a_n)$ such that
$a_i=a_{i+1}+1$ we have
\begin{multline}
\label{dimeq}
\dim H^j(\shn,\Ob(a_1,\ld,a_i,a_{i+1},\ld,a_n))=\\
  =\dim H^j(\shn,\Ob(a_1,\ld,a_{i+1},a_i,\ld,a_n)),\qquad j=0,1,\ld.
\end{multline}
\end{cor}
\begin{proof}
Let $\imath: Y\hk X$ be an embedding, $\E$ a bundle on $X$.
Let $J_Y(\E)$ be a sheaf on $X$:
\begin{equation*}
J_Y(\E)(U)=\{s\in \Gamma(U,\E):\ s|_Y=0\}
\end{equation*}
($U$ is an open set in $X$).
Also denote by $\E^Y$ the following sheaf on $X$:
$\E^Y=\imath_*\imath^* \E$.
We have an exact sequence of the sheaves
$$0\to J_Y(\E)\to\E\to \E^Y\to 0.$$
Note that if $X$ and $Y$ are smooth projective complex algebraic varieties,
$\dim X=n$, $\codim Y=1$, then $J_Y(\E)$ is a locally
free sheaf and thus gives rise to some bundle $\E'$ on $X$, such that
$c_1(\E')=c_1(\E)-[Y]_2$ ($c_1$ is the first Chern class). Here we fix a
notation
$[Y]_2\in H^2(X,\Z)$ for the
image of the class
$[Y]\in H_{2n-2}(X,\Z)$ via the identification
$H_{2n-2}(X,\Z)\simeq H^2(X,\Z)$.
Now let $X=\shn$, $Y=\shin$, $\E=\Ob(a_1,\ld, a_n)$.
We have an exact sequence:
\begin{equation}
\label{es}
0\to J_{\shin}(\Ob(\an))\to\Ob(\an)\to \Ob(\an)^{\shin}\to 0.
\end{equation}
Let us compute the bundle $J_{\shin}(\Ob(\an))$.
Recall that we have fixed a set of generators $[C_i]$, $i=0,\ld,n-1$ of
$H_2(\shn,\Z)$ such that the set of dual generators of
$H_{n-2}(\shn,\Z)$ is $\shin$, $i=1,\ld,n$. It is easy to see that
$c_1(\Ob(a_1,\ld,a_n))([C_j])=a_1+\ld+a_{n-j}$. Note also that
$[\shin]_2=[C_{n-i}]^*$. We obtain
$$c_1(J_{\shin}(\Ob(\an)))=c_1(\Ob(\an))-[C_{n-i}]^*.$$ Thus
$$J_{\shin}(\Ob(\an))=\Ob(a_1,\ld,a_{i-1}, a_i-1,a_{i+1}+1,\ld,a_n).$$
In our case ($a_i=a_{i+1}+1$) we have $a_i-1=a_{i+1}, a_{i+1}+1=a_i$.

From the short exact sequence $(\ref{es})$ we obtain an exact sequence of the
cohomologies:
\begin{multline}
0\to H^0(\Ob(a_1,\ld,a_{i+1},a_i,\ld,a_n))\to H^0(\Ob(\an))\to \\
\to
H^0(\Ob(\an)^{\shin})\to
H^1(\Ob(\ld,a_{i+1},a_i,\ld))\to\\ \to H^1(\Ob(\an))\to
H^1(\Ob(\an)^{\shin})\to\ld
\end{multline}
To prove our lemma it is enough to show that
$H^i(\Ob(\an)^{\shin})=0$ for all ~$i$. Recall that
$\shin$ is fibered over $\sh^{(n-i)}$ with a fiber
$\Pro^1\times \sh^{(i-1)}$. We have proved that the restriction of the bundle
$\Ob(\an)$ to the first factor of the fiber is $\Ob(a_{i+1}-a_i)$. In our
case this
restriction is equal to $\Ob(-1)$. But this bundle has no cohomologies. Thus
the sheaf $\Ob(\an)^{\shin}$ doesn't have any cohomologies too.
Lemma is proved.
\end{proof}

\begin{lem}
\label{sur}
For any $A=(a_1\le \ld\le a_n)\in\N^n$ there exists a surjective map
$\shn\sur \sh_A$.
\end{lem}
\begin{proof}
The idea is to construct such set $B=(b_1 <\ld < b_n)\in\N^n$ that
there exists a permutation $\sigma$ of the factors of $F^{\T (b_n-1)}$
(recall that $F$ is a space of the semi-infinite forms and
$M^A\hk F^{\T (a_n-1)}$) and a vector $v\in F^{\T (b_n-a_n)}$ such that
$\sigma (v_B)=v_A\T v$. If we find such $B$, the existence of the surjective
map $\sh_B\sur \sh_A$ will follow from the definition of the Schubert variety.

Define a map $ar: \N^n\to \N^n$ in the following way. Let $i$
($1\le i< n$) be
a number with a properties $a_i=a_{i+1}$ and $a_j<a_{j+1}$ for $j<i$. Then
$$ar (A)=(a_1,\ld,a_i, a_{i+1}+1,a_{i+2}+1,\ld,a_n+1).$$
Surely, there exist such number $N$ that $ar^N(A)$ has no coinciding
elements. Let $B$ be $ar^N(A)$ with a minimum $N$ with the above property.
From the definition of $v_A$ and $v_B$ one can see that there exists such
vector $v\in F^{\T (b_n-a_n)}$ that $\sigma (v_B)=v_A\T v$ for some permutation
$\sigma$. The lemma is proved.
\end{proof}

\begin{rem}
The above map $\shn\sur \sh_A$
is a resolution of the singularities of the variety $\sh_A$.
\end{rem}

\begin{lem}
\label{ge}
Let $A=(a_1\le a_2\le \ld \le a_n)$, $a_1> 0$. Then
$$\dim H^0(\shn,\Ob(a_1-1,\ld,a_n-1))\ge \Prod_{i=1}^n a_i.$$
\end{lem}
\begin{proof}
First, consider the case $0< a_1 <a_2 <\ld < a_n$. Then
there is an embedding
$\imath_A: \shn\hk \Pro(M^A)$ and $\imath_A^*(\Ob(1))=\Ob(a_1-1,\ld,a_n-1)$.
Thus we have a restriction map
$$H^0(\Pro(M^A),\Ob(1))\to H^0(\shn,\Ob(a_1-1,\ld,a_n-1)).$$
We claim that it is an embedding. In fact, otherwise there exists a linear
homogeneous function $p(x_1,\ld,x_N)$ (where $x_i$ are the coordinates in
$M^A$), which vanishes on $\shn$. Because of the proposition
$(\ref{coordring})$ the latter means that we have some linear condition
on the base elements of $(M^A)^*$.

In spite of the absence of the embedding $\imath_A$ in the case when
there exists $a_i=a_{i+1}$, there is still a map
$\widehat{\imath_A}:\shn\to \Pro(M^A)$, which
is a composition of the surjective map from the lemma $(\ref{sur})$ and
the embedding of $\sh_A$ to $\Pro(M^A)$. This finishes the proof of the lemma.
\end{proof}

\begin{theorem}
\label{maintheorem}
Let $A=(0\le a_1\le\ld\le a_n)$. Then
\begin{equation}
\label{equality}
\dim H^0(\shn, \Ob(\an))=\prod_{i=1}^n (a_i+1)
\end{equation}
and all the higher cohomologies vanish.
\end{theorem}

\begin{proof}
First we will prove the inequality
$$\dim H^0(\shn, \Ob(\an))\le \prod_{i=1}^n (a_i+1)$$ by the induction on $n$.
For $n=1$ it is obvious, because $\sh^{(1)}=\Pro^1$.
Suppose our statement is true for $n-1$. Denote by $M$ the preimage of
some point from $\Pro^1$ via the map $\pi_{n,1}$. Surely,
$M\simeq \sh^{(n-1)}$. Note that $[M]=[\sh^{(n-1)}_n]=[C_0]^*$. Thus we have
the following exact sequence (see the proof of the corollary
$(\ref{coheq})$):
\begin{equation}
\label{exseq}
0\to \Ob(a_1,\ld,a_{n-1}, a_n-1)\to \Ob(\an)\to \Ob(\an)^M\to 0.
\end{equation}
Note also that if $\gamma$ is an embedding of $M$ to $\shn$, then
$$\gamma^* \Ob(a_1,\ld,a_n)=\Ob(a_1,\ld,a_{n-1}).$$ Hence
$H^i(\Ob(\an)^M)\simeq H^i(\Ob(a_1,\ld,a_{n-1}))$. Let us write an exact
sequence
of cohomologies associated with the short exact sequence $(\ref{exseq})$:
\begin{multline}
\label{les}
0\to H^0 (\Ob(a_1,\ld,a_{n-1}, a_n-1))\to H^0(\Ob(\an))\to \\
\to H^0(\Ob(a_1,\ld,a_{n-1}))\to
H^1 (\Ob(a_1,\ld,a_{n-1}, a_n-1))\to \\ \to H^1(\Ob(\an))\to
H^1(\Ob(a_1,\ld,a_{n-1}))\to\\ \to  H^2(\Ob(a_1,\ld,a_{n-1}, a_n-1)) \to\ld
\end{multline}
By the induction assumption we know that
$$\dim H^0(\Ob(a_1,\ld,a_{n-1}))\le\prod_{i=1}^{n-1} (a_i+1).$$
Thus we obtain from the exact sequence $(\ref{les})$ that
\begin{equation}
\label{sum}
\dim H^0(\Ob(\an))\le \prod_{i=1}^{n-1} (a_i+1)+
\dim H^0(\Ob(a_1,\ld, a_{n-1},a_n-1)).
\end{equation}
Now, applying many times the inequality $(\ref{sum})$, using the formula
$(\ref{dimeq})$ for $j=0$ and a fact that $\Ob(0,\ld,0)$ is a trivial bundle
we obtain that
$\dim H^0(\Ob(\an))\le \prod_{i=1}^n (a_i+1)$.
Because of the lemma $(\ref{ge})$ we obtain the equality $(\ref{equality})$.

To clarify the above procedure, let us give an example of the using
of the inequality $(\ref{sum})$ and formula
$(\ref{dimeq})$. Denote by $d(a_1,\ld,a_n)$ the dimension of the space
$H^0(\Ob(\an))$. Let $A=(2,3,4)$. Then we have:
\begin{multline}
\label{ex}
d(2,3,4)\le d(2,3,3)+3\cdot 4\le d(2,3,2)+12+12
=d(2,2,3)+24\le\\
\le d(2,2,2)+9+24 \le d(2,2,1)+9+33=d(2,1,2)+42
=d(1,2,2)+42\le\\
\le d(1,2,1)+6+42=d(1,1,2)+48\le d(1,1,1)+4+48\le d(1,1,0)+4+52=\\
=d(1,0,1)+56=d(0,1,1)+56\le
d(0,1,0)+58=\\
=d(0,0,1)+58\le d(0,0,0)+59=60=3\cdot 4\cdot 5.
\end{multline}

To finish the proof of the theorem we must show that all the higher
cohomologies vanish. Already proved statement $(\ref{equality})$
allows us to rewrite the exact sequence $(\ref{les})$ in the following way
(we use the induction assumption about the vanishing of the higher
cohomologies for $n-1$):
\begin{multline*}
0\to H^1 (\Ob(a_1,\ld,a_{n-1}, a_n-1))\to H^1(\Ob(\an))\to 0\to \\
\to  H^2(\Ob(a_1,\ld,a_{n-1}, a_n-1)) \to H^2(\Ob(\an))\to 0\to \ld.
\end{multline*}
Using this exact sequence and the corollary $(\ref{coheq})$ we obtain the
vanishing of the higher cohomologies. The theorem is proved.
\end{proof}

\begin{cor}
\label{sections}
Let $A=(1\le a_1\le\ld\le a_n)$. Then we have an isomorphism of the
$\slt\T(\C[t]/t^n)$ modules:
$$H^0(\shn,\Ob(a_1-1,\ld,a_n-1))^*\simeq \C^{a_1}*\ld *\C^{a_n}=
M^{(a_1,\ld,a_n)}.$$
\end{cor}
\begin{proof}
Consider the embedding $\imath_{A}:\shn\to \Pro(M^A)$ in the case of the
different $a_i$ and a map $\widehat{\imath_{A}}:\shn\to \Pro(M^A)$ otherwise.
The above theorem states that the map
$$H^0 (\Pro(M^A),\Ob(1))\to H^0(\shn, \Ob(a_1-1,\ld,a_n-1))$$
coming from the equalities
$$\imath_A^* \Ob(1)=\Ob(a_1-1,\ld,a_n-1),\qquad
\widehat{\imath_A}^* \Ob(1)=\Ob(a_1-1,\ld,a_n-1)$$
is an isomorphism.
But $H^0 (\Pro(M^A),\Ob(1))\simeq (M^A)^*$.
\end{proof}

\newcounter{a}

\end{document}